\documentclass[a4paper,12pt]{article}
\usepackage{amssymb}

\newcommand{\R}{\mathbb{R}}

\newcommand{\N}{\mathbb{N}}

\newcommand{\cuad}{{\sqcap\kern-.68em\sqcup}}

\newcommand{\norm}[1]{\|#1\|}

\newtheorem{definition}{Definition}[section]
\newtheorem{theorem}{Theorem}[section]
\newtheorem{proposition}{Proposition}[section]

\newtheorem{lemma}{Lemma}[section]

\newtheorem{remark}{Remark}[section]
\newcommand{\bremark}{\begin{remark} \em}
\newcommand{\eremark}{\end{remark} }

\begin{document}	
\begin{center}
{\Large \bf
Existence, Non-existence, Uniqueness of solutions\\[0.5mm]
 for semilinear  elliptic equations  involving \\[2mm]
  measures concentrated on boundary}

  \bigskip  \medskip

 Huyuan Chen\footnote{hc64@nyu.edu} \quad \quad  Hichem Hajaiej\footnote{hh62@nyu.edu}
\medskip

\begin{abstract}
The purpose of this paper is to study the  weak solutions of the
fractional elliptic problem
\begin{equation}\label{000}
\arraycolsep=1pt
\begin{array}{lll}
 (-\Delta)^\alpha   u+\epsilon g(u)=k\frac{\partial^\alpha\nu}{\partial \vec{n}^\alpha}\quad  &{\rm in}\quad\ \ \bar\Omega,\\[3mm]
 \phantom{(-\Delta)^\alpha   +\epsilon g(u)}
 u=0\quad &{\rm in}\quad\ \ \bar\Omega^c,
 \end{array}
\end{equation}
where $k>0$, $\epsilon=1$ or $-1$,  $(-\Delta)^\alpha$ with  $\alpha\in(0,1)$ is  the fractional Laplacian defined in the principle value sense, $\Omega$ is a bounded $C^2$ open set in $\R^N$ with $N\ge 2$,  $\nu$ is a bounded Radon measure supported in $\partial\Omega$ and $\frac{\partial^\alpha\nu}{\partial \vec{n}^\alpha}$ is defined in the distribution sense, i.e.
$$
\langle\frac{\partial^\alpha\nu}{\partial \vec{n}^\alpha},\zeta\rangle=\int_{\partial\Omega}\frac{\partial^\alpha\zeta(x)}{\partial \vec{n}_x^\alpha}d\nu(x), \qquad \forall\zeta\in C^\alpha(\R^N),
$$
here $\vec{n}_x$ denotes the unit inward normal vector  at $x\in\partial\Omega$.

 In this paper, we prove that (\ref{000}) with $\epsilon=1$ admits a unique weak solution when $g$ is a continuous nondecreasing function satisfying
$$\int_1^\infty (g(s)-g(-s))s^{-1-\frac{N+\alpha}{N-\alpha}}ds<+\infty.$$
Our interest then is to analyse the properties of weak solution  when $\nu=\delta_{x_0}$ with $x_0\in\partial\Omega$, including
the asymptotic behavior near $x_0$ and the limit of weak solutions as $k\to+\infty$.
Furthermore, we show the optimality of the critical value $\frac{N+\alpha}{N-\alpha}$ in a certain sense, by proving the non-existence of weak solutions  when   $g(s)=s^{\frac{N+\alpha}{N-\alpha}}$.

 The final part of this article is  devoted to the  study of  existence for  positive weak solutions
to (\ref{000}) when $\epsilon=-1$ and $\nu$ is a bounded nonnegative Radon measure supported in $\partial\Omega$.
 We employ the Schauder's fixed point theorem to obtain
positive solution under the hypothesis that $g$ is a continuous  function satisfying
 $$\int_1^\infty g(s)s^{-1-\frac{N+\alpha}{N-\alpha}}ds<+\infty.$$
\end{abstract}

\end{center}

 \vspace{1mm}


\bigskip

\setcounter{equation}{0}
\section{Introduction}
\label{sec:intro}

\subsection{Motivation}

In 1991, a fundamental contribution of semilinear elliptic
equations involving measures as boundary data is due to Gmira and V\'{e}ron in \cite{GV}, which studied the weak solutions for
\begin{equation}\label{1.1.1}
\arraycolsep=1pt
\begin{array}{lll}
-\Delta  u+g(u)=0\quad &{\rm in}\quad \Omega,
\\[2mm]\phantom{-----}
u=\mu\quad&{\rm on}\quad \partial\Omega,
\end{array}
\end{equation}
where $\Omega$ is a bounded $C^2$ domain in $\R^N$  and $\mu$ is a bounded Radon measure defined in $\partial\Omega$.
A function $u$ is said to be a weak solution of (\ref{1.1.1})
{\it  if $u\in
L^1(\Omega)$, $g(u)\in L^1(\Omega,\rho_{\partial\Omega} dx)$  and
\begin{equation}\label{1.1.1.0}
\int_\Omega [u(-\Delta)\xi+ g(u)\xi]dx=\int_{\partial\Omega}\frac{\partial\xi(x)}{\partial\vec{n}_x}d\mu(x),\quad \forall\xi\in C^{1.1}_0(\Omega),
\end{equation}
where $\rho_{\partial\Omega}(x)={\rm dist}(x,\partial\Omega)$ and $\vec{n}_x$ denotes the unit inward normal vector  at point $x$.}
Gmira and V´eron  proved that problem (\ref{1.1.1}) admits a unique weak solution  when
 $g$ is a continuous and nondecreasing function satisfying
\begin{equation}\label{14.04}
\int_1^\infty [g(s)-g(-s)]s^{-1-\frac{N+1}{N-1}}ds<+\infty.
\end{equation}
Furthermore, the weak solution of (\ref{1.1.1}) is approached by the classical solutions  of (\ref{1.1.1}) replacing $\mu$ by a
sequence of regular functions $\{\mu_n\}$, which converge to $\mu$ in the distribution sense. Then this subject has been vastly expanded in recent works,
see the papers of Marcus and V\'{e}ron \cite{MV1,MV2,MV3,MV4}, Bidaut-V\'{e}ron  and  Vivier \cite{BV} and reference therein.

A very challenging question consists in  studying the  analogue elliptic  problem involving fractional Laplacian
defined by
$$ (-\Delta)^\alpha  u(x)=\lim_{\varepsilon\to0^+} (-\Delta)_\varepsilon^\alpha u(x),$$
where
$$
(-\Delta)_\varepsilon^\alpha  u(x)=-\int_{\R^N\setminus B_\varepsilon(x)}\frac{ u(z)-
u(x)}{|z-x|^{N+2\alpha}} dz
$$
for $\varepsilon>0$.
The main difficulty comes from how to define the boundary type data.
Given a Radon measure $\mu$ defined in $\partial\Omega$,  it is ill-posed that
$$
\arraycolsep=1pt
\begin{array}{lll}
(-\Delta)^\alpha  u+g(u)=0\quad &{\rm in}\quad \Omega,
\\[2mm]\phantom{------\ }
u=\mu\quad&{\rm on}\quad \partial\Omega,
\\[2mm]\phantom{------\ }
u=0\quad&{\rm in}\quad \bar\Omega^c.
\end{array}
$$
Indeed, let $\{\mu_n\}$ be a sequence of regular functions defined in $\partial\Omega$ converging to the measure $\mu$
and a surprising result is that  there is just zero
solution for
 $$\arraycolsep=1pt
\begin{array}{lll}
(-\Delta)^\alpha  u+g(u)=0\quad &{\rm in}\quad \Omega,
\\[2mm]\phantom{------\ }
u=\mu_n\quad&{\rm on}\quad \partial\Omega,
\\[2mm]\phantom{------\ }
u=0\quad&{\rm in}\quad \bar\Omega^c,
\end{array}
$$
which is in sharp contrast with Laplacian case, where (\ref{1.1.1}) replacing $\mu$ by $\mu_n$ admits a unique nontrivial solution.
On the other hand,
it is also not proper to pose
$$
\arraycolsep=1pt
\begin{array}{lll}
(-\Delta)^\alpha  u+g(u)=0\quad &{\rm in}\quad \Omega,
\\[2mm]\phantom{------\ }
u=\mu\quad&{\rm in}\quad \Omega^c
\end{array}
$$
with $\mu$ being a Radon measure in $\Omega^c$ concentrated on $\partial\Omega$.
In fact,  letting  functions $\{\mu_n\}\subset C^1_0(\Omega^c)$  converging to $\mu$,
the solution $u_n$ of
$$\arraycolsep=1pt
\begin{array}{lll}
(-\Delta)^\alpha  u+g(u)=0\quad &{\rm in}\quad \Omega,
\\[2mm]\phantom{------\ }
u=\mu_n\quad&{\rm in}\quad \Omega^c,
\end{array}$$
is equivalent to the  solution  of
$$\arraycolsep=1pt
\begin{array}{lll}
(-\Delta)^\alpha  u+g(u)=G_{\mu_n}\quad &{\rm in}\quad \Omega,
\\[2mm]\phantom{------\ }
u=0\quad&{\rm in}\quad \Omega^c,
\end{array}$$
where $$G_{\mu_n}(x)=\int_{\Omega^c}\frac{\mu_n(y) }{|x-y|^{N+2\alpha}}dy,\qquad  x\in\Omega,$$
see \cite{CFQ}. It could be seen that
$$\int_\Omega [u_n(-\Delta)^\alpha\xi+ g(u_n)\xi]dx=\int_{\Omega}G_{\mu_n}\xi dx,\quad \forall\xi\in C^2_0(\Omega),$$
Then the limit of $\{u_n\}$ as $n\to\infty$ wouldn't be a weak solution as we desired, similar to (\ref{1.1.1.0}).

Therefore, a totally different point of view has to be found to
 propose the fractional  elliptic  problem involving measure concentrated on boundary.
 Our idea is inspired by the
study of elliptic equations with fractional Laplacian and Radon measure inside of $\Omega$ in \cite{CV1}, where the authors
considered the equations
\begin{equation}\label{1.22}
 \arraycolsep=1pt
\begin{array}{lll}
 (-\Delta)^\alpha  u+h(u)=\nu\quad & \rm{in}\quad\Omega,\\[2mm]
 \phantom{   (-\Delta)^\alpha  +h(u)}
u=0\quad & \rm{in}\quad \Omega^c
\end{array}
\end{equation}
for $\nu\in\mathfrak{M}(\Omega,\rho_{\partial\Omega}^\beta)$ with $ \beta\in[0,\alpha]$  the space of Radon measure $\nu$ in $\Omega$
satisfying $$\int_\Omega\rho_{\partial\Omega}^\beta(x) d|\nu(x)|<+\infty.$$
A function $u$ is said to be a weak solution of (\ref{1.22}), if $u\in
L^1(\Omega)$, $h(u)\in L^1(\Omega,\rho_{\partial\Omega}^\alpha dx)$  and
$$
\int_\Omega [u(-\Delta)^\alpha\xi+ h(u)\xi]dx=\int_{\Omega}\xi(x)d\nu(x),\qquad \forall\xi\in \mathbb{X}_\alpha,
$$
where
$\mathbb{X}_{\alpha}\subset C(\R^N)$ with $\alpha\in(0,1)$  denotes the space of functions
$\xi$ satisfying:\smallskip
\begin{itemize}
\item[]
\begin{enumerate}\item[$(i)$]
${\rm supp}(\xi)\subset\bar\Omega$;
\end{enumerate}
\begin{enumerate}\item[$(ii)$]
 $(-\Delta)^\alpha\xi(x)$ exists for all $x\in \Omega$
and $|(-\Delta)^\alpha\xi(x)|\leq C$ for some $C>0$;
\end{enumerate}
\begin{enumerate}\item[$(iii)$]
there exist $\varphi\in L^1(\Omega,\rho^\alpha_{\partial\Omega} dx)$
and $\varepsilon_0>0$ such that $|(-\Delta)_\varepsilon^\alpha\xi|\le
\varphi$ a.e. in $\Omega$, for all
$\varepsilon\in(0,\varepsilon_0]$.
\end{enumerate}
\end{itemize}
A unique weak
solution of (\ref{1.22}) is obtained when the function $h$ is continuous, nondecreasing and satisfies
$$\int_1^{\infty}(h(s)-h(-s))s^{-1-k_{\alpha,\beta}}ds<+\infty,$$
where
$$
k_{\alpha,\beta}=\left\{
\arraycolsep=1pt
\begin{array}{lll}
\frac{N}{N-2\alpha},\quad &{\rm if}\quad
\beta\in[0,\frac{N-2\alpha}N\alpha],\\[2mm]
\frac{N+\alpha}{N-2\alpha+\beta},\qquad &{\rm if}\quad
\beta\in(\frac{N-2\alpha}N\alpha,\alpha].
\end{array}
\right.
$$

Motivated by the above results, we may approximate $\frac{\partial^\alpha\nu}{\partial \vec{n}^\alpha}$ by a sequence measures defined in $\Omega$
and consider the limit of corresponding weak solutions. To this end, for
a bounded Radon measure defined in $\bar\Omega$ with  support in $\partial\Omega$, we observe that
$$
\langle\frac{\partial^\alpha\nu}{\partial \vec{n}^\alpha},\xi\rangle=\int_{\partial\Omega}\frac{\partial^\alpha\xi(x)}{\partial \vec{n}_x^\alpha}d\nu(x), \quad \xi\in \mathbb{X}_\alpha,
$$
and
$$
\frac{\partial^\alpha\xi(x)}{\partial \vec{n}_x^\alpha}=\lim_{s\to0^+}\frac{\xi(x+s\vec{n}_x)-\xi(x)}{s^\alpha}
=\lim_{s\to0^+}\xi(x+s\vec{n}_x)s^{-\alpha},
$$
so  $\frac{\partial^\alpha\nu}{\partial \vec{n}^\alpha}$ could be approximated by measures $\{t^{-\alpha}\nu_t\}$ with support in $\{x\in\Omega: \rho_{\partial\Omega}(x)=t\}$ generated by $\nu$, see Section 2 for details. Then we consider the limit of weak solutions
as $t\to0^+$ for the problem:
$$
\arraycolsep=1pt
\begin{array}{lll}
(-\Delta)^\alpha  u+g(u)=t^{-\alpha}\nu_t\quad &{\rm in}\quad \Omega,
\\[2mm]\phantom{------\ }
u=0\quad&{\rm in}\quad \Omega^c.
\end{array}
$$
Here  the limit of these weak solutions (if it exists) is called  a weak solution of the following fractional elliptic problem
with  measure concentrated on boundary
$$
\arraycolsep=1pt
\begin{array}{lll}
 (-\Delta)^\alpha   u+g(u)=\frac{\partial^\alpha\nu}{\partial \vec{n}^\alpha}\quad  &{\rm in}\quad\ \ \bar\Omega,\\[3mm]
 \phantom{------\ }
 u=0\quad &{\rm in}\quad\ \ \bar\Omega^c.
 \end{array}
$$
This will be our main focus in this paper.

\subsection{Statement of our problem and main results}

Let $\alpha\in(0,1)$, $g:\R\to\R$ be a continuous function, $\Omega$ be a bounded smooth domain in $\R^N$ with $N\ge 2$ and denote by $\mathfrak{M}^b_{\partial\Omega}(\bar\Omega)$ the bounded Radon measure in $\bar\Omega$ with the support in  $\partial\Omega$.
Our purpose in this article is to investigate the existence, non-existence and uniqueness
of   weak solutions to  semilinear fractional elliptic  problem
\begin{equation}\label{eq 1.1}
\arraycolsep=1pt
\begin{array}{lll}
 (-\Delta)^\alpha   u+\epsilon g(u)=k\frac{\partial^\alpha\nu}{\partial \vec{n}^\alpha}\quad  &{\rm in}\quad\ \ \bar\Omega,\\[3mm]
 \phantom{------\ \  }
 u=0\quad &{\rm in}\quad\ \ \bar\Omega^c,
 \end{array}
\end{equation}
where $\epsilon=1$ or $-1$, $k>0$,   $(-\Delta)^\alpha$  is the fractional Laplacian
and denote $\frac{\partial^\alpha\nu}{\partial \vec{n}^\alpha}$ with $\nu\in\mathfrak{M}^b_{\partial\Omega}(\bar\Omega)$ by
$$
\langle\frac{\partial^\alpha\nu}{\partial \vec{n}^\alpha},\xi\rangle=\int_{\partial\Omega}\frac{\partial^\alpha\xi(x)}{\partial \vec{n}_x^\alpha}d\nu(x), \qquad \xi\in \mathbb{X}_\alpha,
$$
with $\vec{n}_x$ being the unit inward normal vector  at $x$. We call $g$ the absorption nonlinearity if
$\epsilon=1$, otherwise it is called as source nonlinearity.

Before starting our main theorems we make precise the notion of weak solution used in this article.
\begin{definition}\label{weak solution GV}
We say that $u$ is a weak solution of (\ref{eq 1.1}), if $u\in
L^1(\Omega)$, $g(u)\in L^1(\Omega,\rho_{\partial\Omega}^\alpha dx)$  and
$$
\int_\Omega [u(-\Delta)^\alpha\xi+\epsilon g(u)\xi]dx=k\int_{\partial\Omega}\frac{\partial^\alpha\xi(x)}{\partial\vec{n}_x^\alpha}d\nu(x),\qquad \forall\xi\in \mathbb{X}_\alpha.
$$
\end{definition}

We notice that $\mathbb{X}_\alpha\supset C_0^2(\Omega)$  is the test functions space when we study
semilinear fractional elliptic equations involving measures, which plays the same role as $C^{1,1}_0(\Omega)$ for dealing with second order elliptic equations with measures, see \cite{CV1,CV2,CV3,CY}. Moreover, it follows from  \cite[Proposition 1.1]{RS} that $\xi $ is $C^\alpha$ ($\alpha$-H\"{o}lder continuous) in $\R^N$ if $\xi\in \mathbb{X}_\alpha$.

Denote by $G_\alpha$ the Green kernel of $(-\Delta)^\alpha$ in $\Omega\times\Omega$ and by $\mathbb{G}_\alpha[\cdot]$ the
Green operator defined as
$$\mathbb{G}_\alpha[\frac{\partial^\alpha\nu}{\partial\vec{n}^\alpha}](x)=\lim_{t\to0^+}\int_{\partial\Omega} G_\alpha(x,y+t\vec{n}_y)t^{-\alpha}d\nu(y). $$

 Now we are ready to state our first result for problem (\ref{eq 1.1}).
 \begin{theorem}\label{teo 1}
Assume that $\epsilon=1$, $k>0$,  $\nu\in\mathfrak{M}^b_{\partial\Omega}(\bar\Omega)$ and $g$ is a continuous nondecreasing function satisfying
$g(0)\ge0$ and
\begin{equation}\label{g1}
 \int_1^\infty [g(s)-g(-s)]s^{-1-\frac{N+\alpha}{N-\alpha}}ds<+\infty.
\end{equation}
Then

 $(i)$ problem (\ref{eq 1.1}) admits a unique weak solution $u_\nu$;

$(ii)$ the mapping $\nu\to u_\nu$ is increasing and
\begin{equation}\label{1.33}
-k\mathbb{G}_\alpha[\frac{\partial^\alpha\nu_-}{\partial\vec{n}^\alpha}](x)\le u_\nu(x)\le k\mathbb{G}_\alpha[\frac{\partial^\alpha\nu_+}{\partial\vec{n}^\alpha}](x),\qquad x\in\Omega,
\end{equation}
where $\nu_+,\nu_-$ are the positive and negative decomposition of $\nu$ such that $\nu=\nu_+-\nu_-$;

$(iii)$ if we assume additionally that $g$ is $C^{\beta}$ locally in $\R$ with $\beta>0$, then
$u_\nu$ is a classical solution of
\begin{equation}\label{1.1}
 \arraycolsep=1pt
\begin{array}{lll}
 (-\Delta)^\alpha   u+g(u)=0\quad & {\rm in}\quad  \Omega,\\[2mm]
 \phantom{------\ }
u=0\quad & {\rm in}\quad \Omega^c\setminus {\rm supp}(\nu).
\end{array}
\end{equation}

\end{theorem}
\smallskip

\noindent We remark that\\ $(i)$ the second equality in (\ref{1.1}) is understood in the sense that $u=0$ in $\Omega^c\setminus {\rm supp}(\nu)$  and
$u$ is continuous at every point in $\partial\Omega\setminus {\rm supp}(\nu)$;\\ $(ii)$ the uniqueness requires the nondecreasing assumption on nonlinearity $g$,
while the existence also holds without  the nondecreasing assumption on $g$;\\ $(iii)$ (\ref{g1}) is called as integral subcritical condition with
critical value $\frac{N+\alpha}{N-\alpha}$, similar integral subcritical conditions see the references \cite{BV,CV1,CV2,V}.

Applied Theorem \ref{teo 1} when $\nu=\delta_{x_0}$  with $x_0\in\partial\Omega$, problem (\ref{eq 1.1}) admits a unique nonnegative weak solution
when $g$ satisfies the hypotheses in  Theorem \ref{teo 1}.
 Our second goal
is to study the further properties of the  weak solution.

\begin{theorem}\label{teo 2}
 Assume that $\epsilon=1$, $k>0$,   $\nu=\delta_{x_0}$ with $x_0\in\partial\Omega$, $g$ is a  nondecreasing function in $C^\beta$ locally in $\R$
  with $\beta>0$ satisfying $g(0)\ge0$ and (\ref{g1}). Let $u_k$ be the weak solution of
(\ref{eq 1.1}), then

$(i)$
\begin{equation}\label{b k}
\lim_{t\to0^+} \frac{u_k(x_0+t\vec{n}_{x_0})}{ \mathbb{G}_\alpha[\frac{\partial^\alpha\delta_{x_0}}{\partial\vec{n}^\alpha}](x_0+t\vec{n}_{x_0})}=k.
\end{equation}

$(ii)$ if additionally  $g(s)=s^p$ with $p\in(1+\frac{2\alpha}{N},\frac{N+\alpha}{N-\alpha})$,
then the limit of $\{u_k\}$ as $k\to\infty$ exists in $\R^N\setminus\{x_0\}$, denoting $u_\infty$. Moreover, $u_\infty$ is a classical solution of
\begin{equation}\label{1.3}
 \arraycolsep=1pt
\begin{array}{lll}
 (-\Delta)^\alpha   u+u^p=0\quad & {\rm in}\quad  \Omega,\\[2mm]
 \phantom{(-\Delta)^\alpha   +u^p}
u=0\quad & {\rm in}\quad \Omega^c\setminus \{x_0\}
\end{array}
\end{equation}
and satisfies
\begin{equation}\label{b k 01}
c_1 \le u_\infty(x_0+t\vec{n}_{x_0})t^{\frac{2\alpha}{p-1}}\le c_2,\qquad \forall t\in(0,\sigma_0),
\end{equation}
where  $c_2>c_1>0$ and $\sigma_0>0$ small enough.

$(iii)$ if we assume more that $g(s)=s^p$ with $p\in(0,1+\frac{2\alpha}{N}]$, then
$$\lim_{k\to\infty}u_k(x)=+\infty,\qquad \forall x\in \Omega.$$
\end{theorem}

We notice that the limit of  $\{u_k\}$ as $k\to\infty$ blows up every where in $\Omega$ when
$g(u)=u^p$ with $1<p\le1+ \frac{2\alpha}{N}$. This phenomena is different from the Laplacian case, which is caused by the nonlocal
characteristic of the fractional Laplacian.

Theorem \ref{teo 1} and Theorem \ref{teo 2}
show the existence and properties of weak solutions to  (\ref{eq 1.1}) in the subcritical case.
One natural question is what happens in the critical  case, i.e.,
$g(s)=s^p$ with $p\ge \frac{N+\alpha}{N-\alpha}$. The results are given by:

\begin{theorem}\label{teo 4}

 Assume that $\epsilon=1$, $k>0$, $\Omega=B_1(e_N)$ with $e_N=(0,\cdots, 0,1)$,    $\nu=\delta_{0}$  and
$g(s)=s^p$ with $p= \frac{N+\alpha}{N-\alpha}$.
Then  problem (\ref{eq 1.1}) doesn't admit any weak solution.

\end{theorem}

In general, the nonexistence of weak solution is obtained by capacity analysis for second order differential elliptic equations involving
measures, see  \cite{V} and references therein. However, it is a very tough job to attain the nonexistence in the capacity
framework by the nonlocal characteristic and the weak sense of $\frac{\partial^\alpha\delta_0}{\partial \vec{n}^\alpha}$, which is weaker than Radon measure. In the proof of Theorem \ref{teo 4}, we make use of the self-similar property
in the half space.


The last goal of this paper is to  consider the fractional  elliptic problem (\ref{eq 1.1}) with source nonlinearity,
that is, $\epsilon=-1$. In the last decades, semilinear elliptic problems with source nonlinearity and measure data
\begin{equation}\label{eq01}
\arraycolsep=1pt
\begin{array}{lll}
 -\Delta   u=g(u)+k\nu\quad  &{\rm in}\quad\ \ \Omega,\\[3mm]
 \phantom{-\Delta }
 u=\mu\quad &{\rm on}\quad\ \ \partial\Omega,
 \end{array}
\end{equation}
 have attracted numerous interests. There are three basic methods to obtain weak solutions.
The first one is to iterate
$$u_{n+1}=\mathbb{G}_1[g(u_n)]+k\mathbb{G}_1[\nu],\quad \forall n\in\N$$
and look for a function $v$ satisfying
$$v\ge \mathbb{G}_1[g(v)]+k\mathbb{G}_1[\nu].$$
When  $g$ is a pure power source, the existence results could be found in the references \cite{BC,BV,BY,KV,V}.
The second method is to   apply duality argument to derive weak solution  when the mapping $r\mapsto g(r)$ is  nondecreasing,  convex  and continuous, see Baras-Pierre \cite{BP}. These two methods are very difficult to deal with for a general source nonlinearity.
Recently, Chen-Felmer-V\'{e}ron in \cite{CFV} introduced a new method to solve problem (\ref{eq01}) when $g$
is a general nonlinearity, where the authors  employed
Schauder's fixed point theorem to obtain the uniform bound and then to approach the weak solution.

Here we  develop the latter method to attain weak solution of (\ref{eq 1.1}) with $\epsilon=-1$ and the main results state as follows.

\begin{theorem}\label{teo 3}
Let  $\epsilon=-1$, $k>0$ and $\nu\in\mathfrak{M}^b_{\partial\Omega}(\bar\Omega)$ nonnegative with
$\norm{\nu}_{\mathfrak{M}^b(\bar\Omega)}=1$.

$(i)$ Suppose that
\begin{equation}\label{06-08-2}
g(s)\le c_3s^{p_0}+\epsilon,\quad \forall s\ge0,
\end{equation}
for some $p_0\in(0,1]$, $c_3>0$ and $\epsilon>0$.
Assume more that  $c_3$ is small enough when $p_0=1$.

Then problem (\ref{eq 1.1}) admits a nonnegative weak solution $u_\nu$ satisfying
\begin{equation}\label{1.5}
 u_\nu(x)\ge \mathbb{G}_\alpha[\frac{\partial^\alpha\nu}{\partial\vec{n}^\alpha}](x),\qquad \forall x\in\Omega.
\end{equation}

$(ii)$ Suppose that
\begin{equation}\label{1.1+++}
 g(s)\le c_4s^{p_*}+\epsilon,\quad \forall s\in[0,1]
\end{equation}
and
\begin{equation}\label{1.4}
g_\infty:=\int_1^{\infty} g(s)s^{-1-\frac{N+\alpha}{N-\alpha}}ds<+\infty,
\end{equation}
where  $c_4,\epsilon>0$ and $p_*>1$.

Then there exist $k_0,\epsilon_0>0$  depending on $c_4, p_* $ and $ g_\infty$  such that for $k\in[0,k_0)$ and $\epsilon\in(0,\epsilon_0)$, problem
(\ref{eq 1.1}) admits a nonnegative  weak solution $u_\nu$ satisfying (\ref{1.5}).
\end{theorem}

We remark that $(i)$ it does not require any restrictions on parameters $c_3, \epsilon, k$ when $p_0\in(0,1)$ or
on parameters  $ \epsilon, \sigma$ when $p_0=1$;
$(ii)$ the integral subcritical condition (\ref{1.4}) has the same critical value with (\ref{g1}).

The rest of the paper is organized as follows. In Section 2 we study the properties of $\frac{\partial^\alpha \nu}{\partial\vec{n}^\alpha}$.
Section  3 is devoted to prove Theorem \ref{teo 1}. In Section 4  we analyse the properties of the weak solution for problem (\ref{eq 1.1})
when $\nu$ is Dirac mass.
The nonexistence of weak solution in the critical case is addressed  in Section 5. Finally we give the proof of Theorem \ref{teo 3}   in Section 6.
\smallskip
\bigskip


\setcounter{equation}{0}

\section{General measure concentrated on boundary}

In this section, we first build the one-to-one connection between
the Radon measure space $\mathfrak{M}^b_{\partial\Omega}(\bar\Omega)$
and the bounded Radon measure space $\mathfrak{M}^b(\partial\Omega)$.

On the one hand, for any $\mu\in \mathfrak{M}^b(\partial\Omega)$, we denote by $\tilde \mu$ the measure generated by
 $\mu$ extending inside $\Omega$ by zero,
 that is,
 $$\tilde \mu(E):=\mu(E\cap\partial\Omega),\qquad \forall E\subset \bar\Omega\ {\rm Borel\ set}, $$
 then
$\tilde\mu\in \mathfrak{M}^b_{\partial\Omega}(\bar\Omega)$.

On the other hand, let $\tilde \mu\in \mathfrak{M}^b_{\partial\Omega}(\bar\Omega)$,  we see that
$$\tilde \mu(E)=\tilde \mu(E\cap\partial\Omega),\qquad \forall E\subset \bar\Omega\ {\rm Borel\ set}.$$
 Denote by $\mu$ a Radon measure such that
$\mu(F):=\tilde \mu(F),\ F\subset \partial\Omega\ {\rm Borel\ set}$.
Then $\tilde \mu(E)=\mu(E\cap\partial\Omega)$ for any Borel set $E\subset\bar\Omega$
and
$$\norm{\tilde \mu}_{\mathfrak{M}^b(\bar\Omega)}=\norm{\mu}_{\mathfrak{M}^b(\partial\Omega)}.$$

Now we make an approximation of  $\frac{\partial^\alpha\nu}{\partial \vec{n}^\alpha}$
by a sequence Radon measure concentrated on one type of manifolds inside of $\Omega$. Indeed,
we observe that there exists $\sigma_0>0$ small such that
$$\Omega_t:=\{x\in \Omega,\ \rho_{\partial\Omega}(x)>t\}$$ is a $C^2$ domain in $\R^N$
 for $t\in[0,\sigma_0]$ and for any $x\in\partial\Omega_t$, there exists a unique
$x_\partial\in\partial\Omega$ such that
$|x-x_\partial|=\rho_{\partial\Omega}(x)$. Conversely, for any $x\in\partial\Omega$, there exists a unique point $x_t\in\partial\Omega_t$
 such that $|x-x_t|=\rho_{\partial\Omega_t}(x)$, where $t\in(0,\sigma_0)$ and $\rho_{\partial\Omega_t}(x)={\rm dist}(x,\partial\Omega_t)$.
 Then for any Borel set $E\subset \partial\Omega$, there exists unique $E_t\subset \partial\Omega_t$
such that $E_t=\{x_t: x\in E\}$.

In what follows, we always assume that $t\in[0,\sigma_0]$.

Denote by $\nu_t$ a Radon measure generated by $\nu$ as
$$\nu_t(E_t)=\nu(E),$$
and then $\nu_t$ is a bounded Radon measure with support in $\partial\Omega_t$ and
$$\nu_t(E)=\nu_t(E\cap\partial\Omega_t),\qquad \forall E\subset\bar\Omega\  {\rm Borel\ set}.$$
In the distribution sense, we have that
\begin{equation}\label{2.1.2}
\langle\nu_t,f\rangle=\int_{\partial\Omega_t}f(x)d\nu_t(x)=\int_{\partial\Omega}f(x+t\vec{n}_x)d\nu(x),\qquad \forall f\in C_0(\Omega).
\end{equation}
Then  we observe that
\begin{equation}\label{2.1.0}
\{x_t: x\in{\rm supp}(\nu)\}={\rm supp}(\nu_t) \quad{\rm and}\quad\norm{\nu_t}_{\mathfrak{M}^b(\bar\Omega)}=\norm{\nu}_{\mathfrak{M}^b(\bar\Omega)}.
\end{equation}

Now we are able to show an approximation of $\frac{\partial^\alpha\nu}{\partial \vec{n}^\alpha}$.

\begin{proposition}\label{pr 2.1}
The sequence of Radon measures $\{t^{-\alpha}\nu_t\}_t$ converges to $\frac{\partial^\alpha\nu}{\partial \vec{n}^\alpha}$ as
$t\to0^+$ in the following distribution sense:
$$\lim_{t\to0^+}\int_{\partial\Omega_t}\xi(x)t^{-\alpha}d\nu_t(x)=\int_{\partial\Omega}\frac{\partial^\alpha\xi(x)}{\partial \vec{n}^\alpha_x}d\nu(x),\qquad \forall\xi\in \mathbb{X}_\alpha.$$
\end{proposition}
{\bf Proof.} It follows from \cite[Proposition 1.1]{RS},  that $\xi\in C^\alpha(\R^N)$ if $\xi\in \mathbb{X}_\alpha$. This together with the fact that supp$(\xi)\subset\bar\Omega$,
$\frac{\partial^\alpha\xi(x)}{\partial \vec{n}_x^\alpha}$ is well-defined for any $x\in\partial\Omega$
and
for $x_t\in\partial\Omega_t$, implies that there exists a unique $x\in\partial\Omega$ such that
$$x_t=x+t\vec{n}_x\quad{\rm and}\quad |x-x_t|=\rho_{\partial\Omega}(x_t),$$
then
$$\xi(x+t\vec{n}_x)t^{-\alpha}=\frac{\xi(x+t\vec{n}_x)-\xi(x)}{t^\alpha},$$
which implies that
$$\xi(\cdot+t\vec{n})t^{-\alpha}\to \frac{\partial^\alpha\xi(\cdot)}{\partial \vec{n}^\alpha}\quad {\rm as}\quad t\to0^+\quad {\rm in}\quad C(\bar\Omega).$$
Along with (\ref{2.1.2}), we have that
$$ \arraycolsep=1pt
\begin{array}{lll}
|\int_{\partial\Omega_t}\xi(x)t^{-\alpha}d\nu_t(x)-\int_{\partial\Omega}\frac{\partial^\alpha\xi(x)}{\partial \vec{n}_x^\alpha}d\nu(x)|
\\[3mm]\phantom{-----}
=|\int_{\partial\Omega}\xi(x+t\vec{n}_x)t^{-\alpha}d\nu(x)-\int_{\partial\Omega}\frac{\partial^\alpha\xi(x)}{\partial \vec{n}_x^\alpha}d\nu(x)|
\\[3mm]\phantom{-----}
\le\int_{\partial\Omega}|\xi(x+t\vec{n}_x)t^{-\alpha}-\frac{\partial^\alpha\xi(x)}{\partial \vec{n}_x^\alpha}|d|\nu(x)|
\\[3mm]\phantom{-----}\to0\quad{\rm as}\ t\to0^+,
\end{array}
$$
which ends the proof. \qquad$\Box$

\smallskip

We note that Proposition \ref{pr 2.1} shows that $\frac{\partial^\alpha\nu}{\partial \vec{n}^\alpha}$
is approximated  by a sequence Radon measure with support in $\Omega$ in the distribution sense and this provides a new method to derive weak solution
of (\ref{eq 1.1})  by considering the limit of the weak solutions to
$$\arraycolsep=1pt
\begin{array}{lll}
 (-\Delta)^\alpha   u+\epsilon g(u)=kt^{-\alpha}\nu_t\quad & {\rm in}\quad  \Omega,\\[2mm]
\phantom{(-\Delta)^\alpha   +\epsilon g(u)}
 u=0\quad  &{\rm in}\quad  \Omega^c.
\end{array}
$$
To end this section, we give a upper bound for $\mathbb{G}_\alpha[\frac{\partial^\alpha|\nu|}{\partial \vec{n}^\alpha}]$ .
\begin{lemma}\label{lm 2.1}
Let $\nu\in\mathfrak{M}^b_{\partial\Omega}(\bar\Omega)$, then there exists $c_5>0$ such that
 $$\mathbb{G}_\alpha[\frac{\partial^\alpha|\nu|}{\partial \vec{n}^\alpha}](x)\le \int_{\partial\Omega}\frac{c_5}{|x-y|^{N-\alpha}}d|\nu|(y),\qquad x\in\Omega.$$

\end{lemma}
{\bf Proof.} From \cite[Theorem 1.1]{BV}, there exists $c_5>0$ independent of $t$ such that for any $(x,y)\in
\Omega\times\partial\Omega_t$, $x\neq y$,
\begin{equation}\label{annex 01}
G_\alpha(x,y)\le c_5\frac{\rho_{\partial\Omega}^\alpha(y)}{|x-y|^{N-\alpha}}=\frac{c_5 t^{\alpha}}{|x-y|^{N-\alpha}}.
\end{equation}
Then for $x\in\Omega$,
\begin{eqnarray*}
 \mathbb{G}_\alpha[\frac{\partial^\alpha|\nu|}{\partial \vec{n}^\alpha}](x) &=& \lim_{t\to0^+}\int_{\partial\Omega_t} G_\alpha(x,y)t^{-\alpha}d|\nu_t|(y)\\
   &\le &  \lim_{t\to0^+}\int_{\partial\Omega_t} \frac{c_5}{|x-y|^{N-\alpha}}d|\nu_t|(y)
   \\&=&\int_{\partial\Omega} \frac{c_5}{|x-y|^{N-\alpha}}d|\nu|(y).
\end{eqnarray*}
We complete the proof.  \qquad$\Box$

\setcounter{equation}{0}
\section{Absorption Nonlinearity}
\label{sec:existence}

In this section,  our goal is to prove the existence and uniqueness of weak solution for fractional elliptic problem (\ref{eq 1.1}) with $\epsilon=1$.
To this end, we first consider the properties of  weak solution of
\begin{equation}\label{2.0.6}
 \arraycolsep=1pt
\begin{array}{lll}
 (-\Delta)^\alpha   u+g_n(u)=kt^{-\alpha}\nu_t\quad & {\rm in}\quad  \Omega,\\[2mm]
\phantom{(-\Delta)^\alpha   +g_n(u)}
u=0\quad  &{\rm in}\quad  \Omega^c,
\end{array}
\end{equation}
where $t\in(0,\sigma_0)$, $\nu_t$ is given in (\ref{2.1.2})
and  $\{g_n\}$ are a sequence of $C^1$ nondecreasing  functions defined on $\R$
such that $g_n(0)=g(0)\ge 0$,
\begin{equation}\label{06-08-0}
 | g_n|\le g,\quad \sup_{s\in\R}|g_n(s)|=n\quad{\rm and}\quad \lim_{n\to\infty}\norm{g_n-g}_{L^\infty_{loc}(\R)}=0.
\end{equation}

The existence and uniqueness of weak solution to (\ref{2.0.6}) is stated as follows.

\begin{proposition}\label{pr 1}
Assume that $k>0$, $\alpha\in(0,1)$, $g_n$ is a $C^1$ nondecreasing  function satisfying
 $g_n(0)\ge 0$ and (\ref{06-08-0}).
Then for $t\in(0,\sigma_0)$, problem (\ref{2.0.6}) admits a unique weak solution $u_{n,k\nu_t}$ such that
$$|u_{n,k\nu_t}|\le k\mathbb{G}_\alpha[t^{-\alpha}|\nu_t|]\quad{\rm a.e.\ in}\quad \Omega $$
and
\begin{equation}\label{12-08-0}
\norm{g_n(u_{n,k\nu_t})}_{L^1(\Omega,\rho_{\partial\Omega}^\alpha dx)}\le ck\norm{\mathbb{G}_\alpha[t^{-\alpha}|\nu_t|]}_{L^1(\Omega)},
\end{equation}
where $c>0$ independent of $t$, $k$ and $n$.

Furthermore, for any fixed $n\in\N$, $t\in(0,\sigma_0)$ and $k>0$, the mapping $\nu\mapsto u_{n,k\nu_t}$ is increasing.

\end{proposition}
{\bf Proof.} For any $t>0$, we observe that $kt^{-\alpha}\nu_t$ is a bounded Radon measure in $\Omega$ and $g_n$ is bounded, then it follows from
  \cite[Theorem 1.1]{CV1} that problem (\ref{2.0.6}) admits a unique weak solution $u_{n,k\nu_t}$. Moreover, $kt^{-\alpha}\nu_t$
 is increasing with respect to $\nu_t$ and $\nu_t$ is increasing with respect to $\nu$ by the definition of $\nu_t$,  then applying
 \cite[ Theorem 1.1]{CV1}, we have that for any fixed $t\in(0,\sigma_0)$ and $k>0$, the mapping $\nu\mapsto u_{n,k\nu_t}$ is increasing.
\qquad$\Box$
\smallskip

To simplify the notation, we always write $u_{n,k\nu_t}$ by $u_{n,t}$ in this section.
 In order to consider the limit of $\{u_{n,t}\}$ as $t\to0^+$, we introduce some auxiliary lemmas  which are the key steps to obtain  $\{g_n(u_{n,t})\}$ uniformly integrable with respect to $t$.
For $\lambda>0$, let us set
 \begin{equation}\label{Slambda}
  S_\lambda=\{x\in \Omega:\mathbb{G}_{\alpha}[t^{-\alpha}|\nu_t|](x)>\lambda\}\quad{\rm and}\quad m(\lambda)=\int_{S_\lambda} \rho_{\partial\Omega}^\alpha(x)dx.
 \end{equation}

\begin{lemma}\label{lm 0}
For  $\nu\in\mathfrak{M}^b_{\partial\Omega}(\bar\Omega)$ and any $t\in(0, \sigma_0)$,   there exists  $c_6>0$ independent of $t$ such that
\begin{equation}\label{annex -0}
m(\lambda)\le c_6\lambda^{-\frac{N}{N-\alpha}}.
\end{equation}
\end{lemma}
{\bf Proof.}  For $\Lambda>0$ and $y\in\partial\Omega_t$ with $t\in(0, \sigma_0/2)$, we denote
$$A_\Lambda(y)=\{x\in\Omega\setminus\{y\}: G_\alpha(x,y)>\Lambda\}\ \ {\rm
{and}}\quad m_\Lambda(y)=\int_{A_\Lambda(y)}\rho_{\partial\Omega}^\alpha(x) dx.$$
For any $(x,y)\in\Omega\times\partial\Omega_t$, $x\neq y$,
it infers by (\ref{annex 01}) that
\begin{eqnarray*}
A_\Lambda(y)\subset \left\{x\in\Omega\setminus\{y\}: \frac{c_5
t^\alpha}{|x-y|^{N-\alpha}}>\Lambda\right\}\subset B_r(y),
\end{eqnarray*}
where
$r=(\frac{c_5t^\alpha}{\Lambda})^{\frac1{N-\alpha}}$.
Thus, $\rho_{\partial\Omega}(x)\le R_0$ for some $R_0>0$ such that $\Omega\subset B_{R_0}(0)$
and
\begin{equation}\label{annex 1xhw}
m_\Lambda(y)\le R_0^\alpha\int_{B_r(y)}dx\le  c_7t^{\frac{N\alpha}{N-\alpha}}\Lambda^{-\frac{N}{N-\alpha}},
\end{equation}
where $c_7>0$ independent of $t$.


For   $y\in\partial\Omega_t$, we have that
\begin{eqnarray*}
\int_{S_\lambda} G_\alpha(x,y)\rho_{\partial\Omega}^\alpha(x)dx\le
\int_{A_\Lambda(y)}G_\alpha(x,y)\rho_{\partial\Omega}^\alpha(x)dx+\Lambda\int_{S_\lambda}
\rho_{\partial\Omega}^\alpha(x)dx.
\end{eqnarray*}
By integration by parts, we obtain
\begin{eqnarray*}
\int_{A_\Lambda(y)}G_\alpha(x,y)\rho_{\partial\Omega}^\alpha(x)dx=\Lambda m_\Lambda(y)+ \int_\Lambda^\infty m_s(y)ds\le c_8t^{\frac{N\alpha}{N-\alpha}}\Lambda^{1-\frac{N}{N-\alpha}},
\end{eqnarray*}
where $c_8>0$ independent of $t$.
Thus,
\begin{eqnarray*}
\int_{S_\lambda} G_\alpha(x,y)\rho_{\partial\Omega}^\alpha(x)dx\le
c_8t^{\frac{N\alpha}{N-\alpha}}\Lambda^{1-\frac{N}{N-\alpha}}+\Lambda \int_{S_\lambda}
\rho_{\partial\Omega}^\alpha(x)dx.
\end{eqnarray*}
 Choose $\Lambda= t^\alpha(\int_{S_\lambda} \rho_{\partial\Omega}^\alpha(x)dx)^{-\frac{N-\alpha}{N}}$ and then
\begin{eqnarray*}
\int_{S_\lambda} G_\alpha(x,y)\rho_{\partial\Omega}^\alpha(x)dx\le c_9t^\alpha
(\int_{S_\lambda} \rho_{\partial\Omega}^\alpha(x)dx)^{\frac{\alpha}{N}},
\end{eqnarray*}
where $c_9=c_8+1$.
Therefore,
\begin{eqnarray*}
 \int_{S_\lambda}
\mathbb{G}_\alpha [t^{-\alpha}|\nu_t|](x)\rho_{\partial\Omega}^\alpha(x)dx&=&\int_\Omega\int_{S_\lambda}
G_\alpha(x,y)\rho_{\partial\Omega}^\alpha(x)dxt^{-\alpha} d|\nu_t(y)|
\\&\le &c_9\int_\Omega d|\nu_t(y)|(\int_{S_\lambda} \rho_{\partial\Omega}^\alpha(x)dx)^{\frac{\alpha}{N}}
\\&\le& c_9\|\nu\|_{\mathfrak{M}^b(\bar\Omega)} (\int_{S_\lambda} \rho_{\partial\Omega}^\alpha(x)dx)^{\frac{\alpha}{N}}.
\end{eqnarray*}
 As a consequence,
\begin{eqnarray*}
\lambda m(\lambda)\le c_9\|\nu\|_{\mathfrak{M}^b(\bar\Omega)}m(\lambda)^{\frac{\alpha}{N}},
\end{eqnarray*}
which implies (\ref{annex -0}). This ends the proof. \qquad$\Box$

\smallskip
From Lemma \ref{lm 0}, it implies that
\begin{equation}\label{24-08-0}
\norm{\mathbb{G}_\alpha[t^{-\alpha}|\nu_t|]}_{M^{\frac{N}{N-\alpha}}(\Omega,\rho_{\partial\Omega}^\alpha dx)}\le c_9\|\nu\|_{\mathfrak{M}^b(\bar\Omega)},
\end{equation}
where $M^{\frac{N}{N-\alpha}}(\Omega,\rho_{\partial\Omega}^\alpha dx)$ is  Marcinkiewicz space with exponent $\frac{N}{N-\alpha}$. The definition and properties of Marcinkiewicz space see the references \cite{BBC,CC,CV1,V}.

In  next lemma,  the uniformly regularity plays an important role in our approximation of weak solution.
\begin{lemma}\label{lm 1}
Assume that $u_{t}$ is a weak solution
of (\ref{2.0.6}) replacing $g_n$ by $g$, a continuous nondecreasing function satisfying $g(0)\ge 0$.
Then for any compact subsets  $\mathcal{K}\subset\Omega$, there exist $t_0>0$, $\beta>0$ small and $c_{10}>0$ independent of $t$
 such that for $t\in(0,t_0]$,
\begin{equation}\label{2.0.8}
\norm{u_{t}}_{C^{\beta}(\mathcal{K})}\le c_{10} \norm{\nu}_{\mathfrak{M}^b(\bar\Omega)}.
\end{equation}
Moreover, if $g$ is $C^\beta$ locally in $\R$, then there exists $c_{11}>0$ independent of $t$ such that
\begin{equation}\label{2.0.9}
\norm{u_{t}}_{C^{2\alpha+\beta}(\mathcal{K})}\le c_{11} \norm{\nu}_{\mathfrak{M}^b(\bar\Omega)}.
\end{equation}

\end{lemma}
{\bf Proof.} We observe from Proposition \ref{pr 1} that
\begin{equation}\label{3--1}
| u_t|\le \mathbb{G}_\alpha[t^{-\alpha}|\nu_t|]\quad{\rm a.e.\ in}\quad  \Omega.
\end{equation}
For compact set $\mathcal{K}$ in $\Omega$, there exists $t_0>0$ such that
$$\mathcal{K}_{5t_0}\subset \Omega,$$
where $\mathcal{K}_r:=\{x\in\R^N:{\rm dist}(x,\mathcal{K})<r\}$ with $r>0$.
Then $\mathcal{K}_{4t_0}\cap \partial\Omega_t=\O$ for any $t\in(0,t_0]$ and
$$\norm{g(u_t)}_{L^\infty(\mathcal{K}_{3t_0})}\le \norm{g(\mathbb{G}_\alpha[t^{-\alpha}|\nu_t|])}_{L^\infty(\mathcal{K}_{3t_0})}.$$

Since $t^{-\alpha}\nu_t$ is a bounded Radon measure in $\Omega$,  there exists a sequence $\{f_n\}\subset C^2_0(\Omega)$
such that $f_n$ converges to $t^{-\alpha}\nu_t$ in the distribution sense and
 for some  $N_{t_0}>0$ such that for $n\ge N_{t_0}$,
 supp$(f_n)\cap \mathcal{K}_{3t_0}=\O$.

We may assume that $g$ is $C^\beta$ locally in $\R$. (In fact, we can choose a sequence of nondecreasing functions
 $\{g_n\}\subset C^\beta(\R)$ such that $g_n(0)\ge0$, $|g_n(s)|\le |g(s)|$ for $s\in\R$ and $g_n\to g$ locally in $\R$ as $n\to\infty$.)
 Let $w_n$ be the classical solution of
\begin{equation}
 \arraycolsep=1pt
\begin{array}{lll}
 (-\Delta)^\alpha  u+g_n(u)=f_n\quad & {\rm in}\quad\Omega,\\[2mm]
 \phantom{   (-\Delta)^\alpha  + g_n(u)}
u=0\quad & {\rm in}\quad \Omega^c.
\end{array}
\end{equation}
By the uniqueness of weak solution to (\ref{2.0.6}), we obtain that, up to some subsequence,
\begin{equation}\label{12-08.2}
u_t=\lim_{n\to\infty} w_n\quad {\rm a.e.\ in}\quad \Omega.
\end{equation}
 We observe that $0\le w_n= \mathbb{G}_\alpha[f_n]-\mathbb{G}_\alpha[g(w_n)]\le \mathbb{G}_\alpha[f_n]$ and
 $\mathbb{G}_\alpha[f_n]$ converges to $\mathbb{G}_\alpha[t^{-\alpha}|\nu_t|]$ uniformly in any compact set of $\Omega\setminus \partial\Omega_t$ and in $L^1(\Omega)$, then there exists $c_{11}>0$ independent of $n$ and $t$ such that
$$\norm{w_n}_{L^\infty(\mathcal{K}_{3t_0})}\le {c_{11}}\norm{\mathbb{G}_\alpha[t^{-\alpha}|\nu_t|]}_{L^\infty(\mathcal{K}_{3t_0})},\quad  \norm{w_n}_{L^1(\Omega)}\le {c_{11}}\norm{\mathbb{G}_\alpha[t^{-\alpha}|\nu_t|]}_{L^1(\Omega)}.$$
 By  \cite[Lemma 3.1]{CV3}, for $\beta\in(0,2\alpha)$,  there exists $c_{12}>0$ independent of $n$ and $t$, such that
 $$\arraycolsep=1pt
\begin{array}{lll}
\norm{w_n}_{C^{\beta}(\mathcal{K}_{2t_0})} \le c_8[\norm{w_n}_{L^1(\Omega)}+\norm{g(w_n)}_{L^{\infty}(\mathcal{K}_{3t_0})}+\norm{w_n}_{L^{\infty}(\mathcal{K}_{3t_0})}]\\[3mm]
\phantom{} \le c_{12}[\|\mathbb{G}_{\alpha}[t^{-\alpha}|\nu_t|]\|_{L^1(\Omega)}+\norm{\mathbb{G}_{\alpha}[t^{-\alpha}|\nu_t|]}_{L^\infty( \mathcal{K}_{3t_0})}+\|g(\mathbb{G}_{\alpha}[t^{-\alpha}|\nu_t|])\|_{L^\infty( \mathcal{K}_{3t_0})}].
\end{array}
$$
It follows by  \cite[Corollary 2.4]{RS} that there exist $c_{13},c_{14}>0$ such that
\begin{equation}\label{2.0.10.0}
\arraycolsep=1pt
\begin{array}{lll}
\norm{w_n}_{C^{2\alpha+\beta}(\mathcal{K})} \le {c_{13}}[\norm{w_n}_{L^1(\Omega)}+\norm{g(w_n)}_{C^{\beta}(\mathcal{K}_{2t_0})}
+\norm{w_n}_{C^{\beta}(\mathcal{K}_{2t_0})}]\\[3mm]
\phantom{------\ } \le c_{14}[\|\mathbb{G}_{\alpha}[t^{-\alpha}|\nu_t|]\|_{L^1(\Omega)}+\norm{\mathbb{G}_{\alpha}[t^{-\alpha}|\nu_t|]}_{L^\infty( \mathcal{K}_{3t_0})}\\[3mm]
\phantom{-------\ }+\|g\|_{C^\beta([0,\|\mathbb{G}_{\alpha}[t^{-\alpha}|\nu_t|]\|_{L^\infty( \mathcal{K}_{3t_0})}])} \|\mathbb{G}_{\alpha}[t^{-\alpha}|\nu_t|]\|_{C^\beta( \mathcal{K}_{3t_0})}].
\end{array}
\end{equation}
Therefore,  together with (\ref{12-08.2}) and the Arzela-Ascoli Theorem, it follows that
$u_t\in C^{2\alpha+\epsilon}(\mathcal{K})$ for $\epsilon\in(0,\beta)$.
Then $w_n\to u_t$ and $f_n\to 0$ uniformly in any compact subset of $\Omega\setminus\partial\Omega_t$ as $n\to\infty$.
 It infers by   \cite[Lemma 3.1]{CV3} that
$$\arraycolsep=1pt
\begin{array}{lll}
\norm{u_t}_{C^{\beta}(\mathcal{K})} \le c_{8}[\norm{u_t}_{L^1(\Omega)}+\norm{g(u_t)}_{L^{\infty}(\mathcal{K}_{3t_0})}+\norm{u_k}_{L^{\infty}(\mathcal{K}_{3t_0})}]
 \\[2mm]\phantom{-}\le
 c_{12}[\|\mathbb{G}_{\alpha}[t^{-\alpha}|\nu_t|]\|_{L^1(\Omega)}+\norm{\mathbb{G}_{\alpha}[t^{-\alpha}|\nu_t|]}_{L^\infty( \mathcal{K}_{3t_0})}+\|g(\mathbb{G}_{\alpha}[t^{-\alpha}|\nu_t|)]\|_{L^\infty( \mathcal{K}_{3t_0})}].
 \end{array}$$
 We next claim that $\|\mathbb{G}_{\alpha}[t^{-\alpha}|\nu_t|]\|_{L^1(\Omega)}$, $\norm{\mathbb{G}_{\alpha}[t^{-\alpha}|\nu_t|]}_{L^\infty( \mathcal{K}_{3t_0})}$ are uniformly bounded.
In fact,  for $x\in\mathcal{K}$ and $y\in\partial\Omega_t$ with $t\in(0,t_0)$, we have that $|x-y|\ge 3t_0$.
By (\ref{annex 01}), it implies that
\begin{equation}\label{2.1.3}
\arraycolsep=1pt
\begin{array}{lll}
\mathbb{G}_{\alpha}[t^{-\alpha}|\nu_t|](x)\le \int_{\partial\Omega_t}\frac{c_5}{|x-y|^{N-\alpha}}d|\nu_t(y)|
 \\[3mm]\phantom{------} \le c_5t_0^{\alpha-N} \norm{\nu_t}_{\mathfrak{M}^b(\bar\Omega)} = c_5t_0^{\alpha-N} \norm{\nu}_{\mathfrak{M}^b(\bar\Omega)}
 \end{array}
\end{equation}
and
\begin{equation}\label{2.1.4}
\arraycolsep=1pt
\begin{array}{lll}
\norm{\mathbb{G}_{\alpha}[t^{-\alpha}\nu_t]}_{L^1(\Omega)}\le \int_{\Omega}\int_{\partial\Omega_t}\frac{c_5}{|x-y|^{N-\alpha}}d|\nu_t(y)|dx
  \\[3mm]\phantom{------}=\int_{\partial\Omega_t}\int_{\Omega}\frac{c_5}{|x-y|^{N-\alpha}}dxd|\nu_t(y)|
\le c_{15} \norm{\nu}_{\mathfrak{M}^b(\bar\Omega)},
 \end{array}
\end{equation}
which implies that
$$\norm{u_t}_{C^{\beta}(\mathcal{K})}\le c_{15}\norm{\nu}_{\mathfrak{M}^b(\bar\Omega)},$$
where $c_{15}>0$   independent of $t$.

Moreover, if $g$ is $C^\beta$ locally in $\R$, similar to (\ref{2.0.10.0}) it implies by (\ref{2.1.3}) and (\ref{2.1.4})
that
$$\norm{u_t}_{C^{2\alpha+\beta}(\mathcal{K})}\le c_{16}\norm{\nu}_{\mathfrak{M}^b(\bar
\Omega)},$$
where $c_{16}>0$   independent of $t$. We conclude by Theorem 2.2 in \cite{CFQ} that $u_t$ is a classical solution of
\begin{equation}\label{2.2.1}
\arraycolsep=1pt
\begin{array}{lll}
 (-\Delta)^\alpha  u+g(u)=0\quad & {\rm in} \quad\Omega\setminus\partial\Omega_t,\\[2mm]
 \phantom{   (-\Delta)^\alpha  + g(u)}
u=0\quad & {\rm in} \quad \Omega^c.
\end{array}
\end{equation}
This ends the proof.\qquad$\Box$

\begin{proposition}\label{pr 2.01}
Assume that $k>0$
and  $\{g_n\}$ are a sequence of $C^1$ nondecreasing  functions defined on $\R$
such that $g_n(0)=g(0)$ and
(\ref{06-08-0}).
Then problem
\begin{equation}\label{10-08-0}
\arraycolsep=1pt
\begin{array}{lll}
 (-\Delta)^\alpha   u+ g_n(u)=k\frac{\partial^\alpha\nu}{\partial \vec{n}^\alpha}\quad  &{\rm in}\quad\ \ \bar\Omega,\\[3mm]
 \phantom{(-\Delta)^\alpha   + g_n(u)}
 u=0\quad &{\rm in}\quad\ \ \bar\Omega^c
 \end{array}
\end{equation}
admits a unique weak solution $u_{n}$ satisfying
\begin{equation}\label{10-08-1}
-k\mathbb{G}_\alpha[\frac{\partial^\alpha\nu_-}{\partial\vec{n}^\alpha}](x)\le u_n(x)\le k\mathbb{G}_\alpha[\frac{\partial^\alpha\nu_+}{\partial\vec{n}^\alpha}](x),\qquad x\in\Omega,
\end{equation}
where $\nu_+,\nu_-$ are the positive and negative decomposition of $\nu$ such that $\nu=\nu_+-\nu_-$.\\
Furthermore,
\begin{equation}\label{12-08-1}
\norm{g_n(u_{n})}_{L^1(\Omega,\rho_{\partial\Omega}^\alpha dx)}\le k\norm{\mathbb{G}_\alpha[\frac{\partial^\alpha|\nu|}{\partial\vec{n}^\alpha}]}_{L^1(\Omega)}
\end{equation}
and $u_n$ is a classical solution of (\ref{1.1}) replacing  $g$ by $g_n$.
\end{proposition}
{\bf Proof.}
{\it To prove the existence of weak solution.}
Since $\nu_t$ is a bounded Radon measure with supp$(\nu_t)\subset\partial \Omega_t$ for $t\in(0,\sigma_0)$,
then  by Proposition \ref{pr 1}, we have that problem (\ref{2.0.6})
admits a unique weak solution $u_{n,t}$ such that
\begin{equation}\label{24-08-1}
|u_{n,t}|\le \mathbb{G}_\alpha [t^{-\alpha}|\nu_t|]\quad{\rm a.e.\ in}\quad \Omega,\qquad \int_{\Omega}|g_n(u_{n,t})| \rho_{\partial\Omega}^\alpha dx\le k\norm{\mathbb{G}_\alpha [t^{-\alpha}|\nu_t|]}_{L^1(\Omega)}
\end{equation}
and
\begin{equation}\label{2.1.1--}
\int_\Omega [u_{n,t}(-\Delta)^\alpha\xi+g_n(u_{n,t})\xi]dx=\int_{\partial\Omega_t}t^{-\alpha}\xi(x)d\nu_t(x),\quad \forall\xi\in \mathbb{X}_\alpha.
\end{equation}

For any compact set $\mathcal{K} \subset \Omega$, there exists $t_0\in(0,\sigma_0)$ such that
$$\mathcal{K} \subset \Omega_t\quad {\rm and}\quad
{\rm dist}(\mathcal{K},\partial\Omega_{t})\ge t_0,\quad \forall t\in(0,t_0].$$
By Lemma \ref{lm 1}, we observe that   for some $\beta\in(0,\alpha)$
$$\norm{u_{n,t}}_{C^\beta(\mathcal{K})}\le c_5t_0^{-N+2\alpha} \norm{\nu}_{\mathfrak{M}^b(\bar\Omega)}.$$
Therefore, up to some subsequence, there exists $u_n$ such that
$$\lim_{t\to0^+}u_{n,t}=u_n\quad{\rm a.e.\ in}\quad \Omega.$$
Then $g_n(u_{n,t})$ converges to $g_n(u_n)$  almost every in $\Omega$ as $t\to0^+$.
By (\ref{24-08-1}) and (\ref{24-08-0}), we have that $\{u_{n,t}\}_t$ is relatively compact in $L^1(\Omega)$, up to subsequence,
$$u_{n,t}\to u_n\quad {\rm in}\ \ L^1(\Omega)\quad {\rm as}\ t\to0^+$$
and then
$$g_n(u_{n,t})\to g_n(u_n)\quad {\rm in}\ \ L^1(\Omega,\rho^\alpha_{\partial\Omega}dx)\quad {\rm as}\ t\to0^+.$$
By Proposition \ref{pr 2.1},
\begin{eqnarray*}
 \int_{\partial\Omega_t}t^{-\alpha}\xi(x)d\nu_t(x)   \to \int_{\partial\Omega}\frac{\partial^\alpha\xi(x)}{\partial\vec{n}_x^\alpha}d\nu(x)\quad {\rm as}\ t\to0^+,
\end{eqnarray*}
Passing to the limit as
$t\to 0^+$ in the identity (\ref{2.1.1--}),
it implies that
$$\int_\Omega [u_n(-\Delta)^\alpha\xi+g_n(u_n)\xi]dx=k\int_{\partial\Omega}\frac{\partial^\alpha\xi(x)}{\partial\vec{n}_x^\alpha}d\nu(x),\quad \forall\xi\in\mathbb{ X}_\alpha.$$
This implies that $u_n$ is a weak solution of (\ref{10-08-0}).
We see that (\ref{12-08-1}) follows by (\ref{12-08-0}) and Lemma \ref{lm 0}. Moreover,
  by the facts that $u_n=\lim_{t\to0^+} u_{n,t}$ and
$$-k\mathbb{G}_\alpha[t^{-\alpha}\nu_-]\le u_{n,t}\le  k\mathbb{G}_\alpha[t^{-\alpha}\nu_+]\quad {\rm in}\quad \Omega,$$
we have that (\ref{10-08-1}) holds.

\smallskip

{\it To prove that $u_n=0$ in $\Omega^c\setminus{\rm supp}(\nu)$. } Let $x_0\in\partial\Omega\setminus{\rm supp}(\nu)$ and
$x_s=x_0+s\vec{n}_{x_0}$ with $s\in(0,\sigma_0)$. We only have to prove that $\lim_{s\to0^+}u(x_s)=0$.
From \cite[Theorem 1.1]{BV},  for any $(x,y)\in
\Omega\times \partial\Omega_t$, $x\neq y$,
\begin{equation}\label{annex 010}
G_\alpha(x,y)\le c_5\frac{\rho^\alpha_{\partial\Omega}(y)\rho^\alpha_{\partial\Omega}(x)}{|x-y|^N }=c_5\frac{\rho^\alpha_{\partial\Omega}(x)t^{\alpha}}{|x-y|^N }.
\end{equation}
For some $s_0>0$ and any $s\in(0,s_0)$, we observe that ${\rm dist}(x_s,{\rm supp}(\nu))\ge \frac{1}{2}{\rm dist}(x_0,{\rm supp}(\nu))$
and
\begin{eqnarray*}
 \mathbb{G}_\alpha[t^{-\alpha}|\nu_t|](x_s)&\le& c_5\int_{\partial\Omega}\frac{\rho^\alpha_{\partial\Omega}(x_s)}{|x_s-y|^N}d|\nu|(y)
 \\&=& c_5s^\alpha\int_{\partial\Omega\setminus{\rm supp}(\nu)}\frac{1}{|x_s-y|^N}d|\nu|(y)
 \\&\le& c_52^N s^\alpha{\rm dist}(x_0,{\rm supp}(\nu))^{-N}\norm{\nu}_{\mathfrak{M}^b(\bar\Omega)}
 \\&\to&0\quad {\rm as}\quad s\to0^+.
\end{eqnarray*}
Together with the facts that
\begin{equation}\label{facts}
  u_n=\lim_{t\to0^+}u_{n,t}\quad {\rm and}\quad |u_{n,t}|\le \mathbb{G}_\alpha[t^{-\alpha}|\nu_t|],
\end{equation}
we derive that $u_n=0$ in $\Omega^c\setminus{\rm supp}(\nu)$.
\smallskip

{\it To prove the uniqueness of weak solution.}
 Let $u_1,u_2$  be two weak solutions of
(\ref{10-08-0}) and $w=u_1-u_2$. Then $(-\Delta)^\alpha
w=g_n(u_2)-g_n(u_1)$ and $g_n(u_2)-g_n(u_1)\in L^1(\Omega,\rho^\alpha_{\partial\Omega}dx)$.
By Kato's inequatlity, see Proposition 2.4 in \cite{CV1},  for $\xi\in\mathbb{X}_\alpha$, $\xi\ge0$, we have that
\begin{eqnarray*}
\int_\Omega |w|(-\Delta)^\alpha \xi dx+\int_\Omega[g_n(u_1)-g_n(u_2)]{\rm sign}(w)\xi dx\le0.
\end{eqnarray*}
Combining with $\int_\Omega[g_n(u_1)-g_n(u_2)]{\rm sign}(w)\xi dx\ge0$,
then we have
$$w=0\quad {\rm a.e.\ in}\ \ \Omega.$$

{\it Regularity of $u_n$.} Since $g_n$ is $C^1$ in $\R$, then by (\ref{2.0.9}), we have
\begin{equation}\label{2.0.10}
\norm{u_n}_{C^{2\alpha+\beta}(\mathcal{K})}\le c_{17} \norm{\nu}_{\mathfrak{M}^b(\bar\Omega)},
\end{equation}
for any compact set $\mathcal{K}$ and some $\beta\in(0,\alpha)$. Then  $u_n$ is $C^{2\alpha+\beta}$ locally in $\Omega$.
Together with the fact that $u_{n,t}$ is classical solution of (\ref{2.2.1}), we derive by Theorem 2.2 in \cite{CFQ} that
$u_n$ is a classical solution of (\ref{1.1}). \qquad$\Box$
\smallskip

For $\lambda>0$, let us define
 \begin{equation}\label{Slambda}
  \tilde S_\lambda=\{x\in \Omega:\mathbb{G}_{\alpha}[\frac{\partial^\alpha |\nu|}{\partial \vec{n}^\alpha}](x)>\lambda\}\quad{\rm and}\quad \tilde m(\lambda)=\int_{S_\lambda} \rho_{\partial\Omega}^\alpha(x)dx.
 \end{equation}

\begin{lemma}\label{lm 00}
For  $\nu\in\mathfrak{M}^b_{\partial\Omega}(\bar\Omega)$, then there exist $\lambda_0>1$ and $c_{18}>0$ such that for any $\lambda\ge \lambda_0$,
\begin{equation}\label{annex 0}
\tilde m(\lambda)\le c_{18}\lambda^{-\frac{N+\alpha}{N-\alpha}}.
\end{equation}
\end{lemma}
{\bf Proof.} From Lemma \ref{lm 2.1}, we see that
 $$\mathbb{G}_\alpha[\frac{\partial^\alpha|\nu|}{\partial \vec{n}^\alpha}](x)\le \int_{\partial\Omega}\frac{c_5}{|x-y|^{N-\alpha}}d|\nu(y)|,\qquad x\in\Omega.$$
 For $\Lambda>0$ and $y\in\partial\Omega$, we denote
$$\tilde A_\Lambda(y)=\{x\in\Omega: \frac{c_5}{|x-y|^{N-\alpha}}>\Lambda\}\ \ {\rm
{and}}\quad \tilde m_\Lambda(y)=\int_{\tilde A_\Lambda(y)}\rho_{\partial\Omega}^\alpha(x) dx.$$
For any $(x,y)\in\Omega\times\partial\Omega$,
it infers by(\ref{annex 01}) that
\begin{eqnarray*}
\tilde A_\Lambda(y)\subset B_{r_0}(y),
\end{eqnarray*}
where $r_0=(\frac{c_5}{\Lambda})^{\frac1{N-\alpha}}$.

Since $\Omega$ is $C^2$,  there exists $\Lambda_0>1$ such that for $\Lambda>\Lambda_0$ such that
$$\rho_{\partial\Omega}(x)\le |x-y|,\quad \forall x\in \tilde A_\Lambda(y)$$
and
\begin{equation}\label{annex 1xhw}
\tilde m_\Lambda(y)\le \int_{ B_{r_0}(y)}|x-y|^\alpha dx\le  c_{19}\Lambda^{-\frac{N+\alpha}{N-\alpha}}.
\end{equation}


For   $y\in\partial\Omega$, we have that
\begin{eqnarray*}
\int_{\tilde S_\lambda} \frac{c_5}{|x-y|^{N-\alpha}}\rho_{\partial\Omega}^\alpha(x)dx\le
\int_{\tilde A_\Lambda(y)}\frac{c_5}{|x-y|^{N-\alpha}}\rho_{\partial\Omega}^\alpha(x)dx+\Lambda\int_{\tilde S_\lambda}
\rho_{\partial\Omega}^\alpha(x)dx.
\end{eqnarray*}
By integration by parts, we obtain
\begin{eqnarray*}
\int_{\tilde A_\Lambda(y)}\frac{c_5}{|x-y|^{N-\alpha}}\rho_{\partial\Omega}^\alpha(x)dx&=&\Lambda \tilde m_\Lambda(y)+ \int_\Lambda^\infty\tilde m_s(y)ds
\\&\le& c_{20} \Lambda^{1-\frac{N+\alpha}{N-\alpha}},
\end{eqnarray*}
where $c_{20}>0$.
Thus,
\begin{eqnarray*}
\int_{\tilde S_\lambda}\frac{c_5}{|x-y|^{N-\alpha}}\rho_{\partial\Omega}^\alpha(x)dx\le
c_{20}\Lambda^{1-\frac{N+\alpha}{N-\alpha}}+\Lambda \int_{\tilde S_\lambda}
\rho_{\partial\Omega}^\alpha(x)dx.
\end{eqnarray*}
Since $S_{\tilde \lambda_1}\subset S_{\tilde \lambda_2}$ if $\lambda_1\ge \lambda_2$ and
$$\lim_{\lambda\to0^+}\int_{\tilde S_\lambda} \rho_{\partial\Omega}^\alpha(x)dx=0,$$
then there exists $\lambda_0>0$ such that
$$\left(\int_{\tilde S_{\lambda_0}} \rho_{\partial\Omega}^\alpha(x)dx\right)^{-\frac{N-\alpha}{N+\alpha}}\ge\Lambda_0 $$
and for $\lambda\ge \lambda_0$, we may choose $\Lambda= (\int_{\tilde S_\lambda} \rho_{\partial\Omega}^\alpha(x)dx)^{-\frac{N-\alpha}{N+\alpha}}\ge \Lambda_0$ and then
\begin{eqnarray*}
\int_{\tilde S_\lambda} \frac{c_5}{|x-y|^{N-\alpha}}\rho_{\partial\Omega}^\alpha(x)dx\le c_{21}(\int_{S_\lambda} \rho_{\partial\Omega}^\alpha(x)dx)^{\frac{2\alpha}{N+\alpha}},
\end{eqnarray*}
where $c_{21}=c_{20}+1$.
Therefore,
\begin{eqnarray*}
 \int_{\tilde S_\lambda}\mathbb{G}_{\alpha}[\frac{\partial^\alpha |\nu|}{\partial \vec{n}^\alpha}](x)\rho_{\partial\Omega}^\alpha(x)dx&\le &\int_{\partial\Omega}\int_{\tilde S_\lambda}
\frac{c_5}{|x-y|^{N-\alpha}}\rho_{\partial\Omega}^\alpha(x)dx d|\nu(y)|
\\&\le &c_{21}\int_{\partial\Omega} d|\nu(y)|(\int_{\tilde S_\lambda} \rho_{\partial\Omega}^\alpha(x)dx)^{\frac{2\alpha}{N+\alpha}}
\\&\le& c_{21}\|\nu\|_{\mathfrak{M}^b(\bar\Omega)} (\int_{S_\lambda} \rho_{\partial\Omega}^\alpha(x)dx)^{\frac{2\alpha}{N+\alpha}}.
\end{eqnarray*}
 As a consequence,
\begin{eqnarray*}
\lambda \tilde m(\lambda)\le c_{21}\|\nu\|_{\mathfrak{M}^b(\bar\Omega)}\tilde m(\lambda)^{\frac{2\alpha}{N+\alpha}},
\end{eqnarray*}
which implies (\ref{annex 0}). This ends the proof. \qquad$\Box$
\smallskip

To estimate the nonlinearity in $L^1(\Omega,\rho_{\partial\Omega}^\alpha dx)$, we have to introduce an auxiliary lemma as follows.
\begin{lemma}\label{lm 08-09}
Assume that $g:\R_+\mapsto\R_+$ is a continuous function satisfying
\begin{equation}\label{p}
\int_1^{\infty} g(s)s^{-1-p}ds<+\infty
\end{equation}
for some $p>0$.
 Then there is a sequence
real positive numbers $\{T_n\}$ such that
$$\lim_{n\to\infty}T_n=\infty\quad{\rm and}\quad \lim_{n\to\infty}g(T_n)T_n^{-p}=0.$$

Assume additionally that $g$ is nondecreasing, then
$$\lim_{T\to\infty} g(T)T^{-p}=0.$$

\end{lemma}
{\bf Proof.} The first argument see \cite[Lemma 3.1]{CV2} and second see \cite[Lemma 3.1]{CV1}.\quad $\Box$

\smallskip

Now we are ready to prove Theorem \ref{teo 1}.\smallskip

\noindent{\bf Proof of Theorem \ref{teo 1}.}
{\it To prove the existence of weak solution.} Take $\{g_n\}$  a sequence of $C^1$ nondecreasing  functions defined on $\R$
satisfying  $g_n(0)=g(0)$ and
(\ref{06-08-0}).
By Proposition \ref{pr 2.01},  problem (\ref{10-08-0})
admits a unique weak solution $u_{n}$ such that
$$|u_{n}|\le \mathbb{G}_{\alpha}[\frac{\partial^\alpha |\nu|}{\partial \vec{n}^\alpha}]\quad{\rm a.e.\ in}\quad \Omega$$
and
\begin{equation}\label{2.1.1000}
\int_\Omega [u_n(-\Delta)^\alpha\xi+g_n(u_n)\xi]dx=k\int_{\partial\Omega}\frac{\partial^\alpha \xi(x)}{\partial \vec{n}_x^\alpha}d\nu(x),\quad \forall\xi\in \mathbb{X}_\alpha.
\end{equation}

For any compact set $\mathcal{K} \subset \Omega$,
 we observe  from Lemma \ref{lm 1} that   for some $\beta\in(0,\alpha)$,
$$\norm{u_n}_{C^\beta(\mathcal{K})}\le c_{22}  \norm{\nu}_{\mathfrak{M}^b(\bar\Omega)}.$$
Therefore, up to some subsequence, there exists $u_\nu$ such that
$$\lim_{n\to\infty}u_n=u_\nu\quad{\rm a.e.\ in}\ \Omega.$$
Then $ g_n(u_n)$ converge to $g(u_\nu)$ a.e. in $\Omega$ as $n\to\infty$.
By Lemma \ref{lm 00} and (\ref{12-08-1}), we have that
$$u_n\to u_\nu\ {\rm in}\ L^1(\Omega),\quad \norm{g_n(u_n)}_{L^1(\Omega,\rho_{\partial\Omega}^\alpha dx)}\le c_{23}\norm{\mathbb{G}_{\alpha}[\frac{\partial^\alpha |\nu|}{\partial \vec{n}^\alpha}]}_{L^1(\Omega)}$$
and
$$\tilde m(\lambda)\leq c_{18}\lambda^{-\frac{N+\alpha}{N-\alpha}} \ \quad {\rm for}\ \ \  \lambda>\lambda_0,$$
where
$$ \tilde m(\lambda)=\int_{\tilde S_\lambda}\rho_{\partial\Omega}^{\alpha}(x)dx
\quad{\rm  with}\quad \tilde S_\lambda=\{x\in\Omega: \mathbb{G}_{\alpha}[\frac{\partial^\alpha |\nu|}{\partial \vec{n}^\alpha}]>\lambda\}.$$
For any Borel
set $E\subset\Omega$, we have that
$$
\displaystyle\begin{array}{lll}
\displaystyle\int_{E}|g_n(u_n)|\rho_{\partial\Omega}^\alpha(x) dx\le \int_{E\cap\tilde S^c_{\frac{\lambda}{k}}}g\left(k\mathbb{G}_{\alpha}[\frac{\partial^\alpha |\nu|}{\partial \vec{n}^\alpha}]\right)\rho_{\partial\Omega}^\alpha(x) dx+\int_{E\cap \tilde S_{\frac{\lambda}{k}}}g\left(k\mathbb{G}_{\alpha}[\frac{\partial^\alpha |\nu|}{\partial \vec{n}^\alpha}]\right)\rho_{\partial\Omega}^\alpha(x) dx
\\[4mm]\phantom{\int_{E}|g(u_t)|\rho^{\alpha}_{\partial\Omega}(x)dx}
\displaystyle\leq \tilde g\left(\frac{\lambda}{k}\right)\int_E\rho^{\alpha}_{\partial\Omega}(x)dx+\int_{\tilde S_{\frac{\lambda}{k}}}\tilde g\left(k\mathbb{G}_{\alpha}[\frac{\partial^\alpha |\nu|}{\partial \vec{n}^\alpha}]\right)\rho^{\alpha}_{\partial\Omega}(x)dx
\\[4mm]\phantom{\int_{E}|g(u_t)|\rho^{\alpha}_{\partial\Omega}(x)dx}
\displaystyle\leq \tilde
g\left(\frac{\lambda}{k}\right)\int_E\rho^{\alpha}_{\partial\Omega}(x)dx+\tilde m\left(\frac{\lambda}{k}\right) \tilde g\left(\frac{\lambda}{k}\right)+\int_{\frac{\lambda}{k}}^\infty\tilde m(s)d\tilde
g(s),
\end{array}
$$
where $\tilde g(r)=g(|r|)-g(-|r|)$.

On the other hand,
$$\int_{\frac{\lambda}{k}}^\infty \tilde g(s)d\tilde m(s)=\lim_{T\to\infty}\int_{\frac{\lambda}{k}}^T \tilde g(s)d \tilde m(s).
$$
Thus,
$$\displaystyle\begin{array}{lll}
\displaystyle \tilde m\left(\frac{\lambda}{k}\right) \tilde g\left(\frac{\lambda}{k}\right)+ \int_{\frac{\lambda}{k}}^T \tilde m(s)d\tilde g(s) \le c_{24}\tilde g\left(\frac{\lambda}{k}\right)\left(\frac{\lambda}{k}\right)^{-\frac{N+\alpha}{N-\alpha}}+c_{24}\int_{\frac{\lambda}{k}}^T s^{-\frac{N+\alpha}{N-\alpha}}d\tilde g(s)
\\[4mm]\phantom{-----\ \int_{\lambda}^T \tilde g(s)d\omega(s)}\displaystyle
\leq c_{25}T^{-\frac{N+\alpha}{N-\alpha}}\tilde
g(T)+\frac{c_{24}}{\frac{N+\alpha}{N-\alpha}+1}\int_{\frac{\lambda}{k}}^T
s^{-1-\frac{N+\alpha}{N-\alpha}}\tilde g(s)ds.
\end{array}$$
By assumption (\ref{g1}) and Lemma \ref{lm 08-09} with $p=\frac{N+\alpha}{N-\alpha}$,  $T^{-\frac{N+\alpha}{N-\alpha}}\tilde g(T)\to 0$ when $T\to\infty$, therefore,
$$\tilde m\left(\frac{\lambda}{k}\right) \tilde g\left(\frac{\lambda}{k}\right)+ \int_{\frac{\lambda}{k}}^\infty \tilde m(s)\ d\tilde g(s)\leq \frac{c_{24}}{\frac{N+\alpha}{N-\alpha}+1}\int_{\frac{\lambda}{k}}^\infty s^{-1-\frac{N+\alpha}{N-\alpha}}\tilde g(s)ds.
$$
Notice that the above quantity on the right-hand side tends to $0$
when $\lambda\to\infty$. The conclusion follows: for any
$\epsilon>0$ there exists $\lambda>0$ such that
$$\frac{c_{24}}{\frac{N+\alpha}{N-\alpha}+1}\int_{\frac{\lambda}{k}}^\infty s^{-1-\frac{N+\alpha}{N-\alpha}}\tilde g(s)ds\leq \frac{\epsilon}{2}.
$$
For $\lambda$ fixed,  there exists $\delta>0$ such that
$$\int_E\rho_{\partial\Omega}^\alpha(x) dx\leq \delta\Longrightarrow \tilde g\left(\frac{\lambda}{k}\right)\int_E\rho_{\partial\Omega}^\alpha(x) dx\leq\frac{\epsilon}{2},
$$
which implies that $\{g_n\circ u_n\}$ is uniformly integrable in
$L^1(\Omega,\rho_{\partial\Omega}^\alpha dx)$. Then $g_n\circ u_n\to g\circ u_\nu$ in
$L^1(\Omega,\rho_{\partial\Omega}^\alpha dx)$ by Vitali convergence theorem.

Passing to the limit as
$n\to +\infty$ in the identity (\ref{2.1.1000}),
it implies that
$$\int_\Omega [u_\nu(-\Delta)^\alpha\xi+g(u_\nu)\xi]dx=\int_{\partial\Omega}\frac{\partial^\alpha\xi(x)}{\partial\vec{n}^\alpha_x}d\nu(x),\quad \forall\xi\in\mathbb{ X}_\alpha.$$
Then $u_\nu$ is a weak solution of (\ref{eq 1.1}). Moreover,
it follows by the fact
$$-k\mathbb{G}_\alpha[\frac{\partial^\alpha \nu_-}{\partial \vec{n}^\alpha}]\le u_n\le  k\mathbb{G}_\alpha[\frac{\partial^\alpha \nu_+}{\partial \vec{n}^\alpha}]\quad {\rm in}\ \Omega.$$
which, together with $u_\nu=\lim_{n\to+\infty} u_n$, implies  (\ref{1.33}).

 The arguments including $u_n=0$ in $\Omega^c\setminus{\rm supp}(\nu)$,  uniqueness and regularity  follow the proof of Proposition \ref{pr 2.01}. \qquad$\Box$
\smallskip

The proof of the existence of weak solution is divided  into two steps: the first step
is to get weak solution $u_n$ to (\ref{eq 1.1}) with truncated nonlinearity $g_n$ and then to prove the limit of  $\{u_n\}$ as $n\to\infty$
is our desired weak solution. This is due to
the estimate in Lemma \ref{lm 0} where we only could get exponent $\frac{N}{N-\alpha}$ and in the second step, we make use of
Lemma \ref{lm 00}, the critical exponent of the nonlinearity $g$ could be up to $\frac{N+\alpha}{N-\alpha}$.


\setcounter{equation}{0}

\section{Isolated singularity on boundary}
\label{sec:Dirac Mass}
For simplicity, we assume that $x_0=0$ and $\vec{n}_0$ is the unit inward normal vector at the origin  in what follows and   $u_k$ is the weak solution of
(\ref{eq 1.1}).

\subsection{Weak singularity}
In this subsection, we prove Theorem \ref{teo 2} part $(i)$. The regularity refers to  Theorem \ref{teo 1} in the case that $\nu=\delta_{0}$ with $0\in\partial\Omega$ and our main work is to prove (\ref{b k}).
We start our analysis with an auxiliary lemma.

\begin{lemma}\label{lm 3.1}
Under the hypotheses of Theorem \ref{teo 2} part $(i)$, we assume more that $x_s=s\vec{n}_{0}\in\Omega$ for $s>0$ small, then there exists $c_{26}>1$ such that
$$\frac1{c_{26}}s^{-N+\alpha}\le \mathbb{G}_\alpha[\frac{\partial^\alpha \delta_{0}}{\partial \vec{n}^\alpha}](x_s)\le  c_{26}s^{-N+\alpha}$$
and
$$\lim_{s\to0^+}\mathbb{G}_\alpha[g(\mathbb{G}_\alpha[k\frac{\partial^\alpha \delta_{0}}{\partial \vec{n}^\alpha}](x_s))]s^{N-\alpha}=0.$$

\end{lemma}
{\bf Proof.}  It follows by Lemma \ref{lm 2.1} with $\nu=\delta_0$ that
\begin{equation}\label{4.3-1}
\mathbb{G}_\alpha[\frac{\partial^\alpha \delta_0}{\partial \vec{n}^\alpha}](x)\le \frac{c_5}{|x|^{N-\alpha}}, \qquad \forall x\in\Omega,
\end{equation}
in particular,
$$\mathbb{G}_\alpha[\frac{\partial^\alpha \delta_0}{\partial \vec{n}^\alpha}](x_s)\le \frac{c_5}{s^{N-\alpha}}.$$

Let $y_t=t\vec{n}_0$ with $t\in(0,s/2)$, then
$$|y_t-x_s|=s-t> \frac s2=\frac12\max\{\rho_{\partial\Omega}(y_t),\rho_{\partial\Omega}(x_s)\}$$
and apply \cite[Theorem 1.2]{BV} to derive that there exists $c_{27}>0$ such that
\begin{equation}\label{11.04.2}
G_\alpha(x_s,y_t) \ge  c_{27}\frac{\rho^\alpha_{\partial\Omega}(y_t)\rho^\alpha_{\partial\Omega}(x_s)}{|x_s-y_t|^{N}}.
\end{equation}
Thus,
$$\mathbb{G}_\alpha[t^{-\alpha}\delta_{y_t}](x_s)\ge \frac{c_{27} s^{\alpha}}{|s-t|^{N}},$$
which implies that
$$ \mathbb{G}_\alpha[\frac{\partial^\alpha \delta_{0}}{\partial \vec{n}^\alpha}](x_s)\ge \frac{c_{27}}{s^{N-\alpha}}.$$
$(ii)$ By (\ref{4.3-1}) and monotonicity of $g$, we have that
\begin{eqnarray*}
  \mathbb{G}_\alpha[g(\mathbb{G}_\alpha[k\frac{\partial^\alpha \delta_{0}}{\partial \vec{n}^\alpha}])](x_s)s^{N-\alpha} &\le& \int_\Omega G_\alpha(x_s,y)g\left(\frac{c_5k}{|y|^{N-\alpha}}\right)dy s^{N-\alpha} \\
   &\le &  \int_\Omega\frac{c_5}{|x_s-y|^{N-\alpha}} g\left(\frac{c_5k}{|y|^{N-\alpha}}\right)dy s^{N-\alpha} \\
   &=&  c_5 s^{N-\alpha}\left[\int_{ B_{\frac s2}(x_s)}\frac{|y|^\alpha}{|x_s-y|^{N-\alpha}} g\left(\frac{c_5k}{|y|^{N-\alpha}}\right)dy\right.
\\&& \left.+ \int_{ \Omega\setminus B_{\frac s2}(x_s)}\frac{|y|^\alpha}{|x_s-y|^{N-\alpha}} g\left(\frac{c_5k}{|y|^{N-\alpha}}\right)dy\right]
\\&:=&A_1(s)+A_2(s).
\end{eqnarray*}
For $y\in B_{\frac s2}(x_s)$, we have $\frac s2\le |y|\le \frac{3s}2$ and
by applying  Lemma \ref{lm 08-09}, we derive that
\begin{eqnarray*}
A_1(s)&\le&c_5s^{N+\alpha}g\left(\frac{2^{N-\alpha}c_5k}{s^{N-\alpha}}\right)\int_{
B_{1/2}(\vec{n}_0)} \frac{|z|^\alpha}{|\vec{n}_0-z|^{N-\alpha}}dz
\\&=&c_5r^{-\frac{N+\alpha}{N-\alpha}}g\left(2^{N-\alpha}c_5rk\right)\int_{
B_{1/2}(\vec{n}_0)} \frac{|z|^\alpha}{|\vec{n}_0-z|^{N-\alpha}}dz
\\&\to&0\quad{{\rm as}}\ \  r\to+\infty,
\end{eqnarray*}
where $r=s^{\alpha-N}$. We next claim that $A_2(s)\to0$ as $s\to0^+$.
In fact, for $y\in  B_{\frac s2}(0)$,  we see that $|x_s-y|> s/2$ and
\begin{eqnarray*}
 s^{N-\alpha} \int_{  B_{\frac s2}(0)}\frac{|y|^\alpha}{|x_s-y|^{N-\alpha}} g\left(\frac{c_5k}{|y|^{N-\alpha}}\right)dy&\le& 2^{N-\alpha}\int_{ B_{\frac s2}(0)} |y|^\alpha
g\left(\frac{c_5k}{|y|^{N-\alpha}}\right) dy
\\&=&c_{28}\int_0^{\frac s2} r^\alpha g\left(\frac{c_5k}{r^{N-\alpha}}\right)r^{N-1} dr
\\&=&\frac{c_{28}}{N-\alpha}\int_{s^{-\frac1{N-\alpha}}}^\infty \tau^{-1-\frac{N+\alpha}{N-\alpha}} g\left( c_5k\tau\right) d\tau
\\&\to&0\quad{{\rm as}}\ \ s\to0^+,
\end{eqnarray*}
where the converging used (\ref{g1}).
For $y\in  \Omega\setminus \left(B_{\frac s2}(0)\cup B_{\frac s2}(x_s)\right)$, we have that $|y-x_s|>\frac14 |y|$
and
\begin{eqnarray*}
&&s^{N-\alpha} \int_{\Omega\setminus \left(B_{\frac s2}(0)\cup B_{\frac s2}(x_s)\right)}\frac{|y|^\alpha}{|x_s-y|^{N-\alpha}} g\left(\frac{c_5k}{|y|^{N-\alpha}}\right)dy
\\&&\qquad\le s^{N-\alpha}\int_{B_R(0)\setminus B_s(0)} |y|^{2\alpha-N} g\left(\frac{c_5k}{|y|^{N-\alpha}}\right) dy
\\&&\qquad=c_{29}s^{N-\alpha}\int_s^R \tau^{2\alpha-1} g( c_5k\tau^{\alpha-N}) d\tau
\\&&\qquad=c_{29}\frac{s^{2\alpha-1} g(c_5ks^{\alpha-N})}{(N-\alpha)s^{\alpha-N-1}}\qquad \quad{\rm (L'Hospital's\ Rule) }
\\&&\qquad=\frac{c_{29}}{N-\alpha}s^{N+\alpha} g(c_5ks^{\alpha-N})
\\&&\qquad\to0\quad{{\rm as}}\ \ s\to0^+,
\end{eqnarray*}
for some $R>0$ such that $\Omega \subset B_R(0)$ and $c_{29}>0$.
Then
\begin{eqnarray*}
A_2(s)\to 0\quad{{\rm as}}\ \ s\to0^+.
\end{eqnarray*}
Therefore,
\begin{equation}\label{4.4}
\lim_{s\to0^+}\mathbb{G}_\alpha[g(\mathbb{G}_\alpha[k\frac{\partial^\alpha \delta_{0}}{\partial \vec{n}^\alpha}])](x_s)s^{N-\alpha}=0.
\end{equation}
The proof ends. \qquad$\Box$
\smallskip

\noindent{\bf Proof of  Theorem \ref{teo 2}  $(i)$.} The existence, uniqueness and regularity follow by Theorem \ref{teo 1}.
We only need to prove (\ref{b k}) to complete the proof. We observe that
\begin{eqnarray*}
 k\mathbb{G}_\alpha[\frac{\partial^\alpha \delta_{0}}{\partial \vec{n}^\alpha}](x_s) \ge u_k(x_s)
   &\ge&  k\mathbb{G}_\alpha[\frac{\partial^\alpha \delta_{0}}{\partial \vec{n}^\alpha}](x_s)- \mathbb{G}_\alpha[g(u_k)](x_s)
   \\&\ge&  k\mathbb{G}_\alpha[\frac{\partial^\alpha \delta_{0}}{\partial \vec{n}^\alpha}](x_s)-
   \mathbb{G}_\alpha[g(k\mathbb{G}_\alpha[\frac{\partial^\alpha \delta_{0}}{\partial \vec{n}^\alpha}])](x_s),
\end{eqnarray*}
where  $s>0$ small.
Together with Lemma \ref{lm 3.1},  (\ref{b k}) holds.\qquad$\Box$

\subsection{Strong singularity for $p\in(1+\frac{2\alpha}{N},\frac{N+\alpha}{N-\alpha})$}

In this subsection, we consider the limit of $\{u_k\}$ as $k\to\infty$, where $u_k$ is the weak solution of
$$
\arraycolsep=1pt
\begin{array}{lll}
 (-\Delta)^\alpha   u+u^p=k\frac{\partial^\alpha\delta_0}{\partial \vec{n}^\alpha}\quad  &{\rm in}\quad\ \ \bar\Omega,\\[3mm]
 \phantom{-----\ }
 u=0\quad &{\rm in}\quad\ \ \bar\Omega^c,
 \end{array}
 $$
here $0\in\partial\Omega$ and $p\in(1+\frac{2\alpha}{N},\frac{N+\alpha}{N-\alpha})$.
From Theorem \ref{teo 1} $(iii)$, we know that $u_k$ is a classical solution of
\begin{equation}\label{4.101}
\arraycolsep=1pt
\begin{array}{lll}
 (-\Delta)^\alpha   u+u^p=0\quad  &{\rm in}\quad\ \ \Omega,\\[3mm]
 \phantom{-----\ }
 u=0\quad &{\rm in}\quad\ \ \Omega^c\setminus\{0\}.
 \end{array}
\end{equation}
In order to study the limit of $\{u_k\}$ as $k\to\infty$, we have to obtain a super solution of (\ref{4.101}). To this end, we consider the function
\begin{equation}\label{4.2}
  w_p(x)=|x|^{-\frac{2\alpha}{p-1}},\qquad x\in\R^N\setminus\{0\}.
\end{equation}

\begin{lemma}\label{lm 4.1}
Assume that $p\in(1+\frac{2\alpha}{N},\frac{N+\alpha}{N-\alpha})$ and $w_p$ is defined in (\ref{4.2}).
Then there exists $\lambda_0>0$ such that $\lambda_0 w_p$ is a super solution of
(\ref{4.101}).

\end{lemma}
{\bf Proof.} For $p\in(1+\frac{2\alpha}{N},\frac{N+\alpha}{N-\alpha})$, we have that $-\frac{2\alpha}{p-1}\in (-N,-N+2\alpha)$
and  from \cite{FQ}, it shows that there exists $c(p)<0$ such that
$$(-\Delta)^\alpha w_p(x)=c(p)|x|^{-\frac{2\alpha}{p-1}-2\alpha},\quad x\in\R^N\setminus\{0\}, $$
thus, taking $\lambda_0=|c(p)|^{\frac{1}{p-1}}$, we derive that
$$ (-\Delta)^\alpha (\lambda_0 w_p)+ (\lambda_0w_p)^p=0\quad{\rm in}\quad \R^N\setminus\{0\}.$$
Together with $\lambda_0 w_p>0$ in $\Omega^c$, $\lambda_0 w_p$ is a super solution of
(\ref{4.101}).
The proof ends. \qquad$\Box$

\smallskip
We observe that the super solution $\lambda_0w_p$ constructed in Lemma \ref{lm 4.1} could control the asymptotic behavior of $u_\infty$
near the origin, but for $\partial\Omega\setminus\{0\}$,  $\lambda_0w_p$ does not provide enough information for us.
To control the behavior of $u_\infty$ on  $\partial\Omega\setminus\{0\}$, we have to construct new super solutions.
For any given $y_0\in\partial\Omega\setminus\{0\}$, we denote $\eta_0:\R^N\to[0,1]$  a $C^2$ functions such that
\begin{equation}\label{4.5}
\eta_0(x)=\left\{\arraycolsep=1pt
\begin{array}{lll}
  0,\quad  &x\in B_r(y_0),\\[2mm]
 1,\quad &x\in \R^N\setminus B_{2r}(y_0),
 \end{array}
 \right.
\end{equation}
where $r=\frac{|y_0|}{8}$.
\begin{lemma}\label{lm 4.2}
Assume that $p\in(1+\frac{2\alpha}{N},\frac{N+\alpha}{N-\alpha})$ and $w_{\lambda,j}=\lambda\tilde w_p+j \eta_1$,
where $\lambda,j>0$, $\tilde w_p=w_p\eta_0$ in $\R^N$ and $\eta_1=\mathbb{G}_\alpha[1]$.

Then there exist $\lambda_1>0$ and $j_1>0$ depending on $|y_0|$
 such that $w_{\lambda_1,j_1}$ is a super solution of (\ref{4.101}).

\end{lemma}
{\bf Proof.} For $x\in \Omega\setminus B_{4r}(y_0)$, we have that $\tilde w_p(x)=w_p(x)$
and
\begin{eqnarray*}
 (-\Delta)^\alpha \tilde w_p(x)&=& -\lim_{\epsilon\to0^+}\int_{\R^N\setminus B_\epsilon(x)}\frac{\tilde w_p(z)-w_p(x)}{|z-x|^{N+2\alpha}}dz
 \\&=&  (-\Delta)^\alpha  w_p(x)-\lim_{\epsilon\to0^+}\int_{\R^N\setminus B_\epsilon(x)}\frac{\tilde w_p(z)-w_p(z)}{|z-x|^{N+2\alpha}}dz
 \\&\ge&  (-\Delta)^\alpha  w_p(x)-\int_{B_{2r}(y_0)}\frac{w_p(z)}{|z-x|^{N+2\alpha}}dz
 \\&\ge& c(p)|x|^{-\frac{2\alpha}{p-1}-2\alpha} -c_{30}r^{-\frac{2\alpha}{p-1}-2\alpha},
\end{eqnarray*}
where $c_{30}>0$ and the last inequality used the facts $|z-x|\ge 2r$ and $w_p(z)\le r^{-\frac{2\alpha}{p-1}}$.
For  $x\in B_{2r}(0)\setminus\{0\}$, take $\lambda=\lambda_0$ from Lemma \ref{lm 4.1} and $j\ge c_{30}\lambda_0r^{-\frac{2\alpha}{p-1}-2\alpha}$, then we have
\begin{eqnarray*}
 (-\Delta)^\alpha w_{\lambda,j}(x)+ w_{\lambda,j}^p(x)
 &\ge & -c(p)\lambda_0|x|^{-\frac{2\alpha}{p-1}-2\alpha} +w_p^p(x)\ge 0.
\end{eqnarray*}

We observe that there exists $c_{31}>0$ dependent of $r$ such that
$$|(-\Delta)^\alpha \tilde w_p|\le c_{31}\quad{\rm in}\quad  \Omega\setminus B_{2r}(0),$$
then take $j\ge c_{31}\lambda_0$, we have that
$$(-\Delta)^\alpha w_{\lambda_0, j}\ge 0,\quad \forall x\in \Omega\setminus B_{2r}(0). $$
Therefore, letting $\lambda_1=\lambda_0$ and $j_1=\max\{c_{31}\lambda_0, c_{30}\lambda_0r^{-\frac{2\alpha}{p-1}-2\alpha}\}$, we have
that $$ (-\Delta)^\alpha w_{\lambda_1,j_1}+ w_{\lambda_1,j_1}^p\ge 0\quad{\rm in}\quad \Omega.$$
The proof ends. \qquad$\Box$

\small

Let  $x_s=s\vec{n}_0\in\Omega$ and a set
$$A_r=\bigcup_{s\in (0,r)}B_{\frac s8}(x_s).$$
It is obvious that $A_r$ is a cone with the vertex  at the origin.
\begin{lemma}\label{lm 4.3}
Assume that $p\in(0,\frac{N+\alpha}{N-\alpha})$, then there exists $c_{32}>0$ such that for any $x\in A_{r_0}$,
\begin{equation}\label{4.6}
\mathbb{G}_\alpha[(\mathbb{G}_\alpha[\frac{\partial^\alpha \delta_{0}}{\partial \vec{n}^\alpha}])^p](x)
\le \left\{\arraycolsep=1pt
\begin{array}{lll}
c_{32} |x|^{-(N-\alpha)p+2\alpha}\quad  &{\rm if}\quad p\in (\frac{2\alpha}{N-\alpha}, \frac{N+\alpha}{N-\alpha}),\\[1.5mm]
-c_{32} \ln |x| \quad  &{\rm if}\quad p= \frac{2\alpha}{N-\alpha},\\[1.5mm]
c_{32}\quad  &{\rm if}\quad p\in(0,\frac{2\alpha}{N-\alpha}).
\end{array}
 \right.
\end{equation}
\end{lemma}
{\bf Proof.}
Since $\partial\Omega$ is $C^2$, then for $r_0\in(0, 1/2)$ small enough,
we observe that for any $x\in B_{\frac s8}(x_s)$ with $s\in(0,r_0)$,
$$\frac {3s}{4}\le \rho_{\partial\Omega}(x)\le \frac {5s}{4}$$
and for any $t\in(0,\frac s8)$,
$$|x-x_t|\ge \frac{5s}8 \ge \frac12\max\{\rho_{\partial\Omega}(x),\rho_{\partial\Omega}(x_t)\}.$$
Then it follows by \cite[Theorem 1.1, Theorem 1.2]{BV} that there exists $c_{33}>1$ such that
\begin{equation}
\frac1{c_{33}} s^{\alpha-N} t^\alpha\le G_\alpha(x,x_t)\le c_{33}s^{\alpha-N} t^\alpha, \quad \forall x\in B_{\frac s8}(x_s).
\end{equation}
Thus, there exists $c_{34}>0$ independent of $s,t$ such that
$$\frac1{c_{34}}s^{-N+\alpha}\le \mathbb{G}_\alpha[t^{-\alpha}\delta_{x_t}](x)\le  c_{34}s^{-N+\alpha}, \quad \forall x\in B_{\frac s8}(x_s),$$
which implies that
\begin{equation}\label{11.04.1}
 \frac1{c_{34}}s^{-N+\alpha}\le \mathbb{G}_\alpha[\frac{\partial^\alpha \delta_{0}}{\partial \vec{n}^\alpha}](x)\le  c_{34}s^{-N+\alpha},\quad \forall x\in B_{\frac s8}(x_s).
\end{equation}

From Lemma \ref{lm 2.1}, it shows that  for any $x\in\Omega$,
\begin{equation}\label{11.04.3}
 \mathbb{G}_\alpha[\frac{\partial^\alpha \delta_{0}}{\partial \vec{n}^\alpha}](x)\le  c_5 |x|^{-N+\alpha},\qquad \forall x\in \Omega.
\end{equation}
It follows by (\ref{annex 01}) and (\ref{11.04.3}) that
$$ \arraycolsep=1pt
\begin{array}{lll}
\mathbb{G}_\alpha[(\mathbb{G}_\alpha[\frac{\partial^\alpha \delta_{0}}{\partial \vec{n}^\alpha}])^p](x_s)\le c_5^p\int_\Omega G_\alpha(x_s,y)\frac{1}{|y|^{(N-\alpha)p}}dy
 \\[3mm]\phantom{--------\ }  \le c_5^{p+1} \int_\Omega\frac{|y|^\alpha}{|x_s-y|^{N-\alpha}} \frac{1}{|y|^{(N-\alpha)p}}dy
 \\[3mm]\phantom{--------\ }  = c_5^{p+1} s^{2\alpha-(N-\alpha)p}\int_{\tilde \Omega_s}
\frac1{|\vec{n}_0-z|^{N-\alpha}} \frac{1}{|z|^{(N-\alpha)p-\alpha}}dz
 \\[3mm]\phantom{--------\ } = c_5^{p+1}s^{2\alpha-(N-\alpha)p}\left[\int_{\tilde \Omega_s\cap B_{1/2}(\vec{n}_0)}
\frac1{|\vec{n}_0-z|^{N-\alpha}} \frac{1}{|z|^{(N-\alpha)p-\alpha}}dz\right.
 \\[3mm]\phantom{----------------}  \left.+ \int_{\tilde\Omega_s\cap
B_{\frac12}^c(\vec{n}_0)}
\frac1{|\vec{n}_0-z|^{N-\alpha}} \frac{1}{ |z|^{(N-\alpha)p-\alpha}}dz\right]
 \\[3mm]\phantom{--------\ }  :=c_5^{p+1}s^{2\alpha-(N-\alpha)p}[I_1(s)+I_2(s)],
\end{array}$$
where $\Omega_s=\{sz:\ z\in\Omega\}$.

 We observe that
$$ I_1(s)\le c_{35}\int_{ B_{1/2}(\vec{n}_0)}\frac1{|\vec{n}_0-z|^{N-\alpha}} dz\le c_{36}  $$
and since $(N-\alpha)p-\alpha<N$ by $p\in(0,\frac{N+\alpha}{N-\alpha})$, then
\begin{eqnarray*}
I_2(s)&\le& c_{37}\int_{\tilde \Omega_s}\frac{1}{|z|^{(N-\alpha)p-\alpha}(1+|z|)^{N-\alpha}} dz
\\&\le& c_{37}\int_{B_{\frac Rs}(0)\setminus B_{\frac12}(0)}\frac{1}{|z|^{(N-\alpha)p-2\alpha+N}} dz
\\&\le& \left\{\arraycolsep=1pt
\begin{array}{lll}
c_{38} s^{(N-\alpha)p-2\alpha}\quad  &{\rm if}\quad p\in (\frac{2\alpha}{N-\alpha}, \frac{N+\alpha}{N-\alpha}),\\[1.5mm]
-c_{38} \ln s \quad  &{\rm if}\quad p= \frac{2\alpha}{N-\alpha},\\[1.5mm]
c_{38}\quad  &{\rm if}\quad p\in(0,\frac{2\alpha}{N-\alpha}),
\end{array}
 \right.
\end{eqnarray*}
where $c_{35},c_{36},c_{37},c_{38}>0$ and $R>0$ such that $\Omega \subset B_R(0)$. Then (\ref{4.6}) holds. \qquad$\Box$

\medskip

\noindent{\bf Proof of Theorem  \ref{teo 2} $(ii)$.}
For $p\in(1+\frac{2\alpha}{N},\frac{N+\alpha}{N-\alpha})$, we have that
$$-\frac{2\alpha}{p-1}\in(-N,-N+\alpha)$$
and it follows by Lemma \ref{lm 2.1} that
$$u_k(x)\le k\mathbb{G}_\alpha[\frac{\partial^\alpha \delta_0}{\partial \vec{n}^\alpha}](x)\le \frac{c_5k}{|x|^{N-\alpha}},\quad x\in\Omega.$$
Then
$\lim_{x\in\Omega,|x|\to0}\frac{u_k(x)}{ w_p(x)}=0$
and we claim that
$$u_k\le \lambda_0w_p\quad{\rm in}\quad\Omega. $$
In fact, if it fails, then there  exists $z_0\in\Omega$ such that
$$(u_k-\lambda_0w_p)(z_0)=\inf_{\Omega}(u_k-\lambda_0w_p)={\rm ess}\inf_{\R^N}(u_k-\lambda_0w_p)<0.$$
Then we have
$(-\Delta)^\alpha(u_k-\lambda_0w_p)(z_0)<0$, which contradicts
the fact that
$$(-\Delta)^\alpha(u_k-\lambda_0w_p)(z_0)=\lambda_0w_p^p(z_0)-u_k^p(z_0)>0.$$

 By monotonicity of the mapping $k\to u_k$, there holds
$$u_\infty(x):=\lim_{k\to\infty} u_k(x),\quad x\in\R^N\setminus\{0\},$$
 which is a classical solution
of (\ref{4.2}) and
$$u_\infty(x)\le \lambda_0w_p(x)=  \lambda_0 |x|^{-\frac{2\alpha}{p-1}},\quad \forall x\in\Omega.$$
By applying Lemma \ref{lm 4.2}, we obtain that $u_\infty$ is continuous up to the boundary except the origin.

Finally, we claim that there exists $c_{39}>0$ and $t_0<\sigma_0$ such that
\begin{equation}\label{13-08-0}
 u_\infty(x_t)\ge c_{39}t^{-\frac{2\alpha}{p-1}},\quad \forall t\in(0,t_0),
\end{equation}
where $x_t=t\vec{n}_0\in\Omega$.
Indeed, let $r_k=(\sigma^{-1} k)^{\frac{p-1}{(N-\alpha)p-N-\alpha}}$, where $\sigma>0$ will be chosen later, then $k=\sigma r_k^{\frac{(N-\alpha)p-N-\alpha}{p-1}}$
and for $x\in A_{r_0}\cap \left[B_{r_k}(0)\setminus B_{\frac{r_k}{2}}(0)\right]$, we apply Lemma \ref{lm 4.3} with $p\in(1+\frac{2\alpha}{N},\frac{N+\alpha}{N-\alpha})$ that
\begin{eqnarray*}
 u_k(x)&\ge & k\mathbb{G}_\alpha[\frac{\partial^\alpha \delta_{0}}{\partial \vec{n}^\alpha}](x)-k^p\mathbb{G}_\alpha[(\mathbb{G}_\alpha[\frac{\partial^\alpha \delta_{0}}{\partial \vec{n}^\alpha}])^p](x) \\
 &\ge& c_5k|x|^{\alpha-N}[1-c_{40}k^{p-1}|x|^{(\alpha-N)p+\alpha+N}]
  \\&\ge& c_5\sigma r_k^{-\frac{2\alpha }{p-1}}[1-c_{40}\sigma^{p-1} r_k^{p-1}(r_k/2)^{(\alpha-N)p+\alpha+N}]
\\&\ge& c_5\sigma r_k^{-\frac{2\alpha }{p-1}}[1-c_{40}\sigma^{p-1} 2^{(N-\alpha)p-\alpha-N}]
 \\&\ge& \frac{c_5\sigma}{2}|x|^{-\frac{2\alpha }{p-1}},
\end{eqnarray*}
where we choose $\sigma$ such that $c_{40}\sigma^{p-1} 2^{(N-\alpha)p-\alpha-N}=\frac12$.
Then for any $x\in A_{r_0}\cap B_{r_k}^c(0)$, there exists $k>0$ such that
$x\in A_{r_0}\cap [B_{r_k}(0)\setminus B_{\frac{r_k}2}(0)]$ and then
$$u_\infty(x)\ge u_k(x)\ge \frac{c_5\sigma}{2}|x|^{-\frac{2\alpha }{p-1}},\qquad \forall x\in A_{r_0}\cap B_{r_k}^c(0).$$
This ends the proof. \qquad$\Box$

\subsection{The limit of $\{u_k\}$ blows up when $p\in(0,1+\frac{2\alpha}{N}]$}

In this subsection, we derive the blow-up behavior of the limit of $\{u_k\}$  when $p\in(0,1+\frac{2\alpha}{N}]$. To this end,
we first do  precise estimate for $u_k$.
\begin{lemma}\label{lm 3.2}
Assume  that $g(s)=s^p$ with $p\in(1,\frac{N}{N-\alpha}]$ and $u_k$ is the solution of (\ref{eq 1.1}) obtained by Theorem \ref{teo 1}.
Then there exist $c_{41}>0$, $r_0\in(0,\frac14)$ and $\{r_k\}_k\subset(0,r_0)$ satisfying $r_k\to0$ as $k\to\infty$ such that
\begin{equation}\label{4.3.1}
u_k(x)\ge\frac{c_{41}|x|^{-N}}{-\ln(|x|)},\qquad \forall x\in A_{r_0}\cap B_{r_k}^c(0).
\end{equation}
\end{lemma}
{\bf Proof.}
{\it To prove (\ref{4.3.1}) in the case of $p\in (\frac{2\alpha}{N-\alpha},1+\frac{2\alpha}{N})$.} Let $r_j=j^{-\frac1\alpha}$ with $j\in(k_0,k)$, then $j=r_j^{-\alpha}$.
Applying Lemma \ref{lm 4.3} with $p\in (\frac{2\alpha}{N-\alpha},1+\frac{2\alpha}{N})$ and (\ref{11.04.1}), we have that for $x\in A_{r_0}\cap \left[B_{r_j}(0)\setminus B_{\frac{r_j}2}(0)\right]$,
\begin{eqnarray*}
 u_j(x)&\ge & j\mathbb{G}_\alpha[\frac{\partial^\alpha \delta_{0}}{\partial \vec{n}^\alpha}](x)-j^p\mathbb{G}_\alpha[(\mathbb{G}_\alpha[\frac{\partial^\alpha \delta_{0}}{\partial \vec{n}^\alpha}])^p](x) \\
 &\ge& c_{34}^{-1}jr_j^{\alpha-N}-c_{32}j^p|x|^{(\alpha-N)p+2\alpha}
\\&\ge& c_{34}^{-1}r_j^{-N}-c_{32}r_j^{-\alpha p-(N-\alpha)p+2\alpha}
 \\&\ge& \frac{1}{2c_{34}}|x|^{-N},
\end{eqnarray*}
where  the last inequality holds since $-\alpha p-(N-\alpha)p+2\alpha>-N$ and $r_j\to0$ as $j\to\infty$.
Then for any $x\in A_{r_0}\cap B_{r_k}^c(0)$, there exists $j\in (k_0,k)$ such that
$x\in A_{r_0}\cap [B_{r_j}(0)\setminus B_{\frac{r_j}2}(0)]$ and then
$$u_k(x)\ge u_j(x)\ge \frac{1}{2c_{34}}|x|^{-N},\qquad \forall x\in A_{r_0}\cap B_{r_k}^c(0).$$

{\it To prove (\ref{4.3.1}) in the case of $p\in(0,\frac{2\alpha}{N-\alpha}]$.} Let $r_j=j^{-\frac1\alpha}$ with $j\in(k_0,k)$, then $j=r_j^{-\alpha}$
and for $x\in A_{r_0}\cap \left[B_{r_j}(0)\setminus B_{\frac{r_j}2}(0)\right]$, we have that
\begin{eqnarray*}
 u_j(x)&\ge & j\mathbb{G}_\alpha[\frac{\partial^\alpha \delta_{0}}{\partial \vec{n}^\alpha}](x)-j^p\mathbb{G}_\alpha[(\mathbb{G}_\alpha[\frac{\partial^\alpha \delta_{0}}{\partial \vec{n}^\alpha}])^p](x) \\
 &\ge& c_{34}^{-1}j|x|^{\alpha-N}-c_{32}j^p
\\&\ge& c_{34}^{-1}r_j^{-N}-c_{32}r_j^{-\alpha p}
 \\&\ge& \frac{1}{2c_{34}}|x|^{-N},
\end{eqnarray*}
where  the last inequality holds since $-\alpha p>-N$ and $r_j\to0$ as $j\to\infty$.
For any $x\in A_{r_0}\cap B_{r_k}^c(0)$, there exists $j\in (k_0,k)$ such that
$x\in A_{r_0}\cap [B_{r_j}(0)\setminus B_{\frac{r_j}2}(0)]$ and then
$$u_k(x)\ge u_j(x)\ge \frac{1}{2c_{34}}|x|^{-N},\qquad \forall x\in A_{r_0}\cap B_{r_k}^c(0).$$

{\it To prove (\ref{4.3.1}) in the case of $p=1+\frac{2\alpha}{N}$.} Let $\rho_j=j^{-\frac1\alpha}$ and $r_j=\frac{\rho_j}{[-\log(\rho_j)]^{\frac1\alpha}}$, then $j=\rho_j^{-\alpha}$
and applied Lemma \ref{lm 4.3} for $x\in A_{r_0}\cap \left[B_{r_j}(0)\setminus B_{\frac{r_j}{2}}(0)\right]$,
\begin{eqnarray*}
 u_j(x)&\ge & j\mathbb{G}_\alpha[\frac{\partial^\alpha \delta_{0}}{\partial \vec{n}^\alpha}](x)-j^p\mathbb{G}_\alpha[(\mathbb{G}_\alpha[\frac{\partial^\alpha \delta_{0}}{\partial \vec{n}^\alpha}])^p](x) \\
 &\ge& c_{34}^{-1}j|x|^{\alpha-N}-c_{32}j^p|x|^{(\alpha-N)p+2\alpha}
\\&\ge& c_{34}^{-1}\rho_j^{-N} (-\log \rho_j)^{\frac{N-\alpha}{\alpha}}-c_{42} \rho_j^{-N} (-\log \rho_j)^{\frac{(N-\alpha)p-2\alpha}{\alpha}}
 \\&=& c_{34}^{-1}\rho_j^{-N} (-\log \rho_j)^{\frac{N-\alpha}{\alpha}}\left[1-c_{42}(-\log \rho_j)^{\frac{(N-\alpha)p-2\alpha}{\alpha}-\frac{N-\alpha}{\alpha}}\right]
 \\&\ge& c_{34}^{-1}\frac{r_j^{-N}}{-\log \rho_j}[1-c_{42}(-\log \rho_j)^{\frac{(N-\alpha)p-N-\alpha}{\alpha}}]
 \\&\ge& \frac{c_{34}|x|^{-N}}{-2\log |x|},
\end{eqnarray*}
where $c_{42}>0$ and we used the facts that $\log(\rho_j)\le c\log r_j\le c\log |x|$ and $\frac{(N-\alpha)p-N-\alpha}{\alpha}<0$.
Then for any $x\in A_{r_0}\cap B_{r_k}^c(0)$, there exists $j\in (k_0,k)$ such that
$x\in A_{r_0}\cap [B_{r_j}(0)\setminus B_{\frac{r_j}2}(0)]$ and then
$$u_k(x)\ge \frac{c_{34}|x|^{-N}}{-2\log |x|},\qquad x\in A_{r_0}\cap B_{r_k}^c(0).$$
The proof ends. \qquad$\Box$
\smallskip

\noindent{\bf Proof of Theorem  \ref{teo 2} $(iii)$.}
It derives by Lemma \ref{lm 3.2} that
\begin{equation}\label{3.2.3}
\pi_k:=\int_{B_{r_0}(0)}u_k(x)\ge c_{41}\int_{A_{r_0}\cap B_{r_k}^c(0)}\frac{|x|^{-N}}{-\log|x|}dx\to\infty\quad {\rm as}\ k\to\infty.
\end{equation}
Fix $y_0\in \Omega\setminus  \bar B_{r_0}(0)$, it follows by Lemma 2.4 in \cite{CY} that problem
\begin{equation}\label{5.1}
\arraycolsep=1pt
\begin{array}{lll}
 (-\Delta)^\alpha  u+u^p=0 \quad & {\rm in}\quad  B_{\varrho_0}(y_0),\\[2mm]
 \phantom{  (-\Delta)^\alpha  +u^p}
u=0  \quad & {\rm in}\quad \R^N \setminus (B_{\varrho_0}(y_0)\cup B_{r_0}(0)),\\[2mm]
\phantom{ (-\Delta)^\alpha  +u^p}
u=u_k  \quad & {\rm in}\quad B_{r_0}(0)
\end{array}
\end{equation}
admits a unique solution $w_k$, where $\varrho_0=\min\{\rho_{\partial\Omega}(y_0),|y_0|-r_0\}$.
By Lemma 2.2 in \cite{CY},
\begin{equation}\label{4.1.3}
 u_{k}\ge w_k\quad {\rm in}\quad B_{\varrho_0}(y_0).
\end{equation}
Let  $\tilde w_k=w_k-u_k\chi_{B_{r_0}(0)},$
then $\tilde w_k=w_k$ in $B_{\varrho_0}(y_0)$ and for $x\in B_{\varrho_0}(y_0)$,
$$
\arraycolsep=1pt
\begin{array}{lll}
(-\Delta)^\alpha \tilde w_k(x) =
-\lim_{\epsilon\to0^+}\int_{B_{\varrho_0}(y_0)\setminus B_\epsilon(x)}\frac{w_k(z)-w_k(x)}{|z-x|^{N+2\alpha}}dz
\\[3mm]\phantom{--------}+\lim_{\epsilon\to0^+}\int_{B_{\varrho_0}^c(y_0)\setminus B_\epsilon(x)}\frac{w_k(x)}{|z-x|^{N+2\alpha}}dz
\\[3mm]\phantom{------}
=-\lim_{\epsilon\to0^+}\int_{\R^N\setminus B_\epsilon(x)}\frac{w_k(z)-w_k(x)}{|z-x|^{N+2\alpha}}dz +\int_{B_{r_0}(0)}\frac{u_k(z)}{|z-x|^{N+2\alpha}}dz
\\[3mm]\phantom{------}\ge(-\Delta)^\alpha w_k(x)+c_{42}\pi_k,
\end{array}
$$
where $c_{42}=(|y_0|+r_0)^{-N-2\alpha}$ and the last inequality follows by the fact of $$|z-x|\le |x|+|z|\le |y_0|+r_0\quad
{\rm for}\ z\in B_{\frac14}(0),\ x\in B_{\frac14}(y_0).$$
Therefore,
\begin{eqnarray*}
(-\Delta)^\alpha \tilde w_k(x)+\tilde w_k^p(x) &\ge&  (-\Delta)^\alpha w_k(x)+w_k^p(x)+ c_{42}\pi_k \\
        &=&c_{42}\pi_k, \qquad x\in B_{\varrho_0}(y_0),
     \end{eqnarray*}
that is,  $\tilde w_k$  is a super solution of
\begin{equation}\label{4.1.2}
\arraycolsep=1pt
\begin{array}{lll}
\displaystyle (-\Delta)^\alpha  u+u^p=c_{42}\pi_k \quad & {\rm in}\quad  B_{\varrho_0}(y_0),\\[2mm]
 \phantom{  (-\Delta)^\alpha  +u^{p,}}
u=0  \quad & {\rm in}\quad B_{\varrho_0}^c(y_0).
\end{array}
\end{equation}
Let $\eta_1$ be the solution of
$$
\arraycolsep=1pt
\begin{array}{lll}
 (-\Delta)^\alpha  u=1 \quad & {\rm in}\quad  B_{\varrho_0}(y_0),\\[2mm]
 \phantom{  (-\Delta)^\alpha  }
u=0  \quad & {\rm in}\quad  B^c_{\varrho_0}(y_0).
\end{array}
$$
Then
 $(c_{42}\pi_k)^{\frac1p} \frac{\eta_1}{2\max_{\R^N}\eta_1}$ is sub solution of (\ref{4.1.2}) for $k $ large enough. By Lemma 2.2 in \cite{CY},  we have that
 $$\tilde w_k(x)\ge (c_{42}\pi_k)^{\frac1p} \frac{\eta_1(x)}{2\max_{\R^N}\eta_1},\quad \forall  x\in B_{\varrho_0}(y_0),$$
 which implies that
 $$w_k(y)\ge c_{43} (c_{42}\pi_k)^{\frac1p},\qquad \forall y\in B_{\frac{\varrho_0}{2}}(y_0),$$
 where $c_{43}=\min_{x\in B_{\varrho_0}(y_0)}\frac{\eta_1(x)}{2\max_{\R^N}\eta_1}$.
 Therefore, (\ref{4.1.3}) and (\ref{3.2.3}) imply that $$\lim_{k\to\infty}u_{k}(y)\ge \lim_{k\to\infty}w_k(y)=\infty,\qquad\forall y\in B_{\frac{\varrho_0}{2}}(y_0).$$
Similarly, we can prove
$$\lim_{k\to\infty}u_{k}(y)\ge \lim_{k\to\infty}w_k(y)=\infty,\qquad\forall y\in \Omega.$$
The proof ends.\qquad$\Box$

\setcounter{equation}{0}
\section{Nonexistence in the critical case}


In this section, we prove the nonexistence in the critical case. To this end, we consider the weak solution to  elliptic problem
\begin{equation}\label{14-08-0}
\arraycolsep=1pt
\begin{array}{lll}
 (-\Delta)^\alpha   u+u^{\frac{N+\alpha}{N-\alpha}}=k\frac{\partial^\alpha\delta_0}{\partial e_N^\alpha}\quad  &{\rm in}\quad\ \ \overline{\R^N_+},\\[3mm]
 \phantom{(-\Delta)^\alpha  +u^p}
 u=0\quad &{\rm in}\quad\ \ \R^N_-,
 \end{array}
\end{equation}
where  $\R^N_+=\R^{N-1}\times\R_+$ and $e_N=(0,\cdots,0,1)$.

\begin{definition}\label{weak definition}
A function $u\in L^1(\R^N,\mu dx)$ is a weak solution of (\ref{14-08-0}) if $u^p\in L^1(\R^N,\rho^\alpha \mu dx)$   and
\begin{equation}\label{weak sense}
\int_{\R^N_+} [u(-\Delta)^\alpha\xi+u^{\frac{N+\alpha}{N-\alpha}}\xi]dx=\frac{\partial^\alpha \xi(0)}{\partial e_N^\alpha},\quad \forall\xi\in \mathbb{X}_{\alpha,\R^N_+},
\end{equation}
where  $\mu(x)=\frac1{1+|x|^{N+2\alpha}}$, $\rho(x)=\min\{1,\rho_{\partial\Omega}(x)\}$ and $\mathbb{X}_{\alpha,\R^N_+}\subset C(\R^N)$ is the space of functions
$\xi$ satisfying:\smallskip

\noindent (i) the support of $\xi$ is a compact set in $\bar\R^N_+$;\smallskip

\noindent(ii) $(-\Delta)^\alpha\xi(x)$ exists for any $x\in\R^N_+$ and there exists $c>0$ such that
 $$|(-\Delta)^\alpha\xi(x)|\leq  c\mu(x),\quad \forall x\in\R^N_+;$$

\noindent(iii) there exist $\varphi\in L^1(\R^N_+,\rho^\alpha dx)$
and $\varepsilon_0>0$ such that $|(-\Delta)_\varepsilon^\alpha\xi|\le
\varphi$ a.e. in $\R^N_+$, for all
$\varepsilon\in(0,\varepsilon_0]$.\smallskip
\end{definition}

Let $\mathbb{G}_{\alpha,\R^N_+}$ the Green's function on $\R^N_+\times\R^N_+$ and
\begin{equation}\label{5.2}
\Gamma_\alpha(x)=\lim_{t\to0}t^{-\alpha}\mathbb{G}_{\alpha,\R^N_+}(x,te_N).
\end{equation}

\begin{lemma}\label{lm 5.1}
Let $\Gamma_\alpha$ defined in (\ref{5.2}), then
\begin{equation}\label{5.1.1}
\arraycolsep=1pt
\begin{array}{lll}
(-\Delta)^\alpha \Gamma_\alpha=\frac{\partial^\alpha \delta_0}{\partial e_N^\alpha}\qquad &{\rm in}\quad\bar \R^N_+,\\[2mm]
 \phantom{--- }
\Gamma_\alpha=0  \quad & {\rm in} \quad \R^N_-.
\end{array}
\end{equation}
Moreover,
\begin{equation}\label{5.3}
\Gamma_\alpha(x)=|x|^{-N+\alpha}\Gamma_\alpha\left(\frac{x}{|x|}\right),\qquad x\in\R^N,
\end{equation}
and
$$
\Gamma_\alpha\left(\frac{x}{|x|}\right)\ \left\{\arraycolsep=1pt
\begin{array}{lll}
>0 \quad & {\rm if}\quad  x\in \R^N_+,\\[2mm]
=0  \quad & {\rm if}\quad x\not\in \R^N_+.
\end{array}
\right.
$$
\end{lemma}
{\bf Proof.} We observe that
$$(-\Delta)^\alpha_x t^{-\alpha}\mathbb{G}_{\alpha,\R^N_+}(x,te_N)=t^{-\alpha}\delta_{te_N}$$
and
$$\lim_{t\to0^+}\langle t^{-\alpha}\delta_{te_N},\xi\rangle=\frac{\partial^\alpha \xi(0)}{\partial e_N^\alpha},\qquad \forall \xi\in\mathbb{X}_{\alpha,\R^N_+}.$$
Then (\ref{5.1.1}) holds in the weak sense. By the regularity results,  $\Gamma_\alpha$ is a solution of
\begin{equation}\label{5.1.2}
\arraycolsep=1pt
\begin{array}{lll}
(-\Delta)^\alpha \Gamma_\alpha=0\qquad &{\rm in}\quad\R^N_+,\\[2mm]
 \phantom{--- }
\Gamma_\alpha=0  \quad & {\rm in} \quad\overline{\R^N_-}\setminus\{0\}.
\end{array}
\end{equation}
Let $\Gamma_{\alpha,\lambda}(x)=\lambda^{N-\alpha}\Gamma_{\alpha}(\lambda x)$
 and  $\xi_\lambda(x)=\xi(x/\lambda)$ for $\xi\in\mathbb{X}_{\alpha,\R^N_+}$, then we have that
\begin{eqnarray*}
\int_{\R^N_+}\Gamma_{\alpha,\lambda}(-\Delta)^\alpha \xi dx&=&\lambda^\alpha\int_{\R^N_+}\Gamma_{\alpha}(z)(-\Delta)^\alpha \xi_\lambda(x)dx, \\
   &=&  \lambda^\alpha \frac{\partial^\alpha \xi_\lambda(0)}{\partial e_N^\alpha},
\end{eqnarray*}
which implies that
$$\int_{\R^N_+}\Gamma_{\alpha,\lambda}(-\Delta)^\alpha \xi dx=\frac{\partial^\alpha \xi(0)}{\partial e_N^\alpha}.$$
By the uniqueness, we derive that
$$\lambda^{N-\alpha}\Gamma_{\alpha}(\lambda x)=\Gamma_{\alpha}(x),$$
which, choosing $\lambda=\frac1{|x|}$, implies (\ref{5.3}).
The last argument is obvious.
\qquad$\Box$

\begin{theorem}\label{teo 4.1}
Let  $k>0$,  then problem (\ref{14-08-0})
 has no any weak solution.

\end{theorem}
{\bf Proof.} If there exists a weak solution $u_k$ to (\ref{14-08-0}), then we observe that
$$u_k>0\qquad {\rm in}\quad \R^N_+.$$
By Maximum Principle, we have that
\begin{equation}\label{5.5}
  u_k\le k\Gamma_\alpha\qquad {\rm in}\quad \R^N.
\end{equation}
Denoting
$$u_\infty=\lim_{k\to\infty}u_k\qquad {\rm in}\quad \R^N.$$
We claim that
\begin{equation}\label{3.2.1}
u_\infty(x)=|x|^{\alpha-N}u_\infty(\frac{x}{|x|}),\quad \forall x\in\R^N\setminus\{0\}.
\end{equation}
Indeed, let
$$\tilde u_\lambda(x)=\lambda^{N-\alpha}u_k(\lambda x),\quad \forall x\in \R^N\setminus\{0\}. $$
By direct computation, we have that for $x\in\R^N_+$,
\begin{eqnarray}
 (-\Delta)^\alpha  \tilde u_\lambda(x) +\tilde u_\lambda^{\frac{N+\alpha}{N-\alpha}}(x)
  &=&\lambda^{N+\alpha}[(-\Delta)^\alpha u_k(\lambda x)  +  u_k^{\frac{N+\alpha}{N-\alpha}}(\lambda x)] \nonumber \\
  &=&0.\label{13-09-5}
\end{eqnarray}
Moreover, for $f\in C_0^1(\R^N_+)$,
\begin{eqnarray*}
 \langle(-\Delta)^\alpha  \tilde u_\lambda +\tilde u_\lambda^{\frac{N+\alpha}{N-\alpha}}, f\rangle&=&\lambda^{N+\alpha}\int_{\R^N}
 [(-\Delta)^\alpha u_k(\lambda x)  +  u_k^{\frac{N+\alpha}{N-\alpha}}(\lambda x)]f(x)dx\nonumber
  \\&=&\lambda^{ \alpha }\int_{\R^N}  [(-\Delta)^\alpha u_k(z)  +  u_k^{\frac{N+\alpha}{N-\alpha}}(z)]f\left(\frac{z}{\lambda}\right)dz\nonumber
   \\  &=&\lambda^{\alpha}k\frac{\partial^\alpha f(0)}{\partial e_N^\alpha}.
\end{eqnarray*}
Thus,
\begin{equation}\label{3.2.2}
  (-\Delta)^\alpha  \tilde u_\lambda +\tilde u_\lambda^{\frac{N+\alpha}{N-\alpha}}=\lambda^{\alpha}k\frac{\partial^\alpha \delta_0}{\partial e_N^\alpha}\quad {\rm in}\ \ \R^N_+.
\end{equation}
We observe that $\lim_{|x|\to\infty}\tilde u_\lambda(x)=0$
and  $u_{k\lambda^{\alpha }}$ is the unique weak solution of
(\ref{14-08-0}) with $k$ replaced by $\lambda^\alpha k$,
then  for $x\in \R^N\setminus\{0\}$,
\begin{equation}\label{21-10-1}
u_{k\lambda^{\alpha}}(x)=\tilde u_\lambda (x)=\lambda^{N-\alpha}u_k(\lambda x)
\end{equation}
and letting $k\to\infty$ we have that
$$
u_{\infty}(x)=\lambda^{N-\alpha}u_\infty(\lambda x),\qquad \forall x\in \R^N\setminus\{0\},
$$
which implies (\ref{3.2.1}) by taking $\lambda=|x|^{-1}$.

Combine (\ref{5.3}), (\ref{5.5}) and (\ref{21-10-1}), then we have that
 \begin{eqnarray*}
  u_{k\lambda^{\alpha}}(x) \le  \lambda^{N-\alpha}k\Gamma_\alpha(\lambda x) = k\Gamma_\alpha(x),\quad \forall x\in \R^N.
 \end{eqnarray*}
Thus,
$$u_\infty(x)\le k\Gamma_\alpha(x),\quad \quad \forall x\in \R^N.$$
 By arbitrary of $k$, it implies that
$$u_\infty\equiv0,$$
then $u_1\equiv 0$ in $\R^N$,
which is impossible.
\qquad$\Box$ \medskip

\noindent{\bf Proof of Theorem \ref{teo 4}.} Without loss generality, we let $k=1$, $0\in\partial\Omega$ and $e_N$ is the unit normal vector pointing inside of $\Omega$ at $0$. If
$$
\arraycolsep=1pt
\begin{array}{lll}
 (-\Delta)^\alpha   u+u^{\frac{N+\alpha}{N-\alpha}}=\frac{\partial^\alpha\delta_0}{\partial e_N^\alpha}\quad  &{\rm in}\quad\ \ \bar\Omega,\\[3mm]
 \phantom{(-\Delta)^\alpha  +u^{\frac{N+\alpha}{N-\alpha}}}
 u=0\quad &{\rm in}\quad\ \ \bar\Omega^c
 \end{array}
$$
admits a solution weak $v_1$, we claim that there is a weak solution of (\ref{14-08-0}),
then the contradiction is obtained from Theorem \ref{teo 4.1}.

In fact, we may assume that
$$\Omega=B_1(e_N)\quad{\rm and}\quad B_m=B_m(me_N).$$ Then
$$\Omega\subset B_m\subset B_{m+1}\quad {\rm and}\quad \lim_{m\to\infty} B_m=\R^N_+.$$
Let
$$v_m(x)=m^{\alpha-N}v_1(\frac xm),\quad x\in\R^N.$$
By direct computation, $v_m$ is a weak solution of
\begin{equation}\label{15-08-0-0}
\arraycolsep=1pt
\begin{array}{lll}
 (-\Delta)^\alpha   u+u^{\frac{N+\alpha}{N-\alpha}}=\frac{\partial^\alpha\delta_0}{\partial e_N^\alpha}\quad  &{\rm in}\quad\ \ \bar B_m,\\[3mm]
 \phantom{(-\Delta)^\alpha  +u^{\frac{N+\alpha}{N-\alpha}}}
 u=0\quad &{\rm in}\quad\ \ \bar B_m^c,
 \end{array}
\end{equation}

\smallskip

\noindent{\it  We next show that $v_m\le v_{m+1}$ in $\R^N$. }
From Proposition \ref{pr 1},
\begin{equation}\label{08-14-00}
\arraycolsep=1pt
\begin{array}{lll}
 (-\Delta)^\alpha   u+u^{\frac{N+\alpha}{N-\alpha}}=t^{-\alpha}\delta_{te_N}\quad  &{\rm in}\quad\ \  B_m,\\[3mm]
 \phantom{(-\Delta)^\alpha  +u^{\frac{N+\alpha}{N-\alpha}}}
 u=0\quad &{\rm in}\quad\ \  B_m^c
 \end{array}
\end{equation}
admits a unique weak solution, denoting  $v_{m,t}$.
Choose a sequence nonnegative functions $\{f_{m,i}\}_{i\in\N}\subset C^1(\R^N)$ with support $B_1(e_N)$
such that $f_{m,i}\rightharpoonup t^{-\alpha}\delta_{te_N}$ as $i\to\infty$ in the distribution sense. Let $v_{m,i,t}$   be  the unique solution
of
\begin{equation}\label{08-14-00-0}
\arraycolsep=1pt
\begin{array}{lll}
 (-\Delta)^\alpha   u+u^{\frac{N+\alpha}{N-\alpha}}=f_{m,i}\quad  &{\rm in}\quad\ \  B_m,\\[3mm]
 \phantom{(-\Delta)^\alpha  +u^{\frac{N+\alpha}{N-\alpha}}}
 u=0\quad &{\rm in}\quad\ \  B_m^c
 \end{array}
\end{equation}
and by Maximum Principle, see \cite[Lemma 2.3]{CY}, derive that
$$
v_{m,i,t}\le \tilde v_{m+1,i,t}\quad{\rm in}\quad \R^N.
$$
Together with the facts that
 $v_{m,i,t}\to v_{m,t}$ a.e. in $\R^N$ and $v_{m+1,i,t}\to v_{m+1,t}$ a.e. in $\R^N$ as $i\to\infty$, we obtain that
\begin{equation}\label{3.2}
v_{1,t}\le v_{m,t}\le  v_{m+1,t}\quad{\rm a.e.\ in}\ \ \R^N
\end{equation}
and
$$\int_{B_m}v_{m,t}^{\frac{N+\alpha}{N-\alpha}}\rho^\alpha dx< \norm{\mathbb{G}_{\alpha, B_m}[f_{m,i}]}_{L^1(\Omega,\ \rho^\alpha dx)},$$
which implies that
\begin{equation}\label{08-14-0-1}
\arraycolsep=1pt
\begin{array}{lll}
 (-\Delta)^\alpha  u+u^{\frac{N+\alpha}{N-\alpha}}=\frac{\partial^\alpha \delta_0}{\partial e_N^\alpha}\quad  &{\rm in}\quad\ \  \bar B_m,\\[3mm]
 \phantom{(-\Delta)^\alpha  +u^{\frac{N+\alpha}{N-\alpha}}}
 u=0\quad &{\rm in}\quad\ \  \bar B_m^c
 \end{array}
\end{equation}
admits a solution $v_m$  for any $m\in\N$
and
\begin{equation}\label{3.3}
 v_{m}\le v_{m+1}\quad {\rm a.e.\ in}\ \ \R^N.
\end{equation}
We observe that
\begin{equation}\label{3.4}
0\le v_{m}\le \mathbb{G}_{\alpha, B_{m}}[\frac{\partial^\alpha \delta_0}{\partial e_N^\alpha}]\le \frac{c_5}{|x|^{N-\alpha}} \quad {\rm a.e.\ in}\ \ \R^N
\end{equation}
and
$$\int_{B_m}v_{m}^{\frac{N+\alpha}{N-\alpha}}\rho^\alpha dx< \norm{\mathbb{G}_{\alpha, B_m}[\frac{\partial^\alpha \delta_0}{\partial e_N^\alpha}]}_{L^1(B_m,\ \rho^\alpha dx)}.$$
By (\ref{3.3}) and (\ref{3.4}), we see that the limit of  $\{v_m\}$ exists, denoted it by $w_1$. Hence,
\begin{equation}\label{2.1.1}
0\le w_1\le \mathbb{G}_{\alpha,\R^N_+}[\frac{\partial^\alpha \delta_0}{\partial e_N^\alpha}]\quad {\rm a.e.\ in}\ \ \R^N
\end{equation}
and
$$\int_{\R^N_+}w_1^{\frac{N+\alpha}{N-\alpha}}\rho^\alpha dx< \norm{\mathbb{G}_{\alpha,\R^N_+}[\frac{\partial^\alpha \delta_0}{\partial e_N^\alpha}]}_{L^1(\R^N_+,\rho^\alpha\mu dx)},$$
which implies that $w_1\in L^1(\R^N,\ \mu dx)$.
Thus, $v_m\to w_1$  in $L^1(\R^N,\ \rho^\alpha \mu dx)$ as $m\to\infty$.

For $\xi\in \mathbb{X}_{\alpha,\R^N_+}$, there exists $N_0>0$ such that for any $m\ge N_0$,
 $${\rm supp}(\xi)\subset \bar B_m,$$
which implies that
$\xi\in \mathbb{X}_{\alpha,B_m}$ and then
\begin{equation}\label{3.6}
\int_{\R^N_+} [v_m(-\Delta)^\alpha\xi+v_m^{\frac{N+\alpha}{N-\alpha}}\xi]dx= \frac{\partial^\alpha \xi(0)}{\partial e_N^\alpha}.
\end{equation}
By \cite[Lemma 3.1 ]{CY},
$$|(-\Delta)^\alpha\xi(x)|\le \frac{c_{9}\norm{\xi}_{L^\infty(\Omega)}}{1+|x|^{N+2\alpha}},\quad \forall x\in \R^N_+.$$
Thus,
\begin{equation}\label{3.7}
 \lim_{m\to\infty}\int_{\R^N_+} v_m(x)(-\Delta)^\alpha\xi(x) dx=\int_{\R^N_+} w_1(x)(-\Delta)^\alpha\xi(x) dx.
\end{equation}

By (\ref{2.1.1}) and increasing monotonicity of $v_m$,
for any $n\ge N_0$,
\begin{equation}\label{3.8}
 \lim_{m\to\infty}\int_{\R^N_+} v_m^{\frac{N+\alpha}{N-\alpha}}\xi(x) dx=\int_{\R^N_+} w_1^{\frac{N+\alpha}{N-\alpha}}\xi(x) dx.
\end{equation}
Combining (\ref{3.7}), (\ref{3.8}) and taking $m\to\infty$ in (\ref{3.6}), we obtain that
\begin{equation}\label{3.10}
\int_{\R^N_+} \left[w_1(-\Delta)^\alpha\xi+w_1^{\frac{N+\alpha}{N-\alpha}}\xi\right]dx=\frac{\partial^\alpha \xi(0)}{\partial e_N^\alpha}.
\end{equation}
Since $\xi\in \mathbb{X}_{\alpha,\R^N_+}$ is arbitrary,  $w_1$
is a weak solution of (\ref{14-08-0}).\qquad$\Box$


\setcounter{equation}{0}
\section{Forcing nonlinearity}
This section is devoted to consider problem  (\ref{eq 1.1}) when $\epsilon=-1$, we call it as forcing case.
In order to derive the existence of weak solution to  (\ref{eq 1.1}) with forcing nonlinearity, we first introduce the following  propositions.

\begin{proposition}\label{general} \cite[Proposition 2.2]{CFV}
Let $\alpha\in(0,1]$, $\beta\in[0,\alpha]$ and $\nu\in\mathfrak{M}(\Omega,\rho^\beta_{\partial\Omega})$, then there exists $c_{44}>0$ such that
\begin{equation}\label{annex 00}
\|\mathbb{G}_\alpha[\nu]\|_{M^{p_\beta^*}(\Omega,\rho^\beta_{\partial\Omega}
dx)}\le c_{44}\|\nu\|_{\mathfrak{M}(\Omega,\rho^\beta_{\partial\Omega})},
\end{equation}
where $p_\beta^*=\frac{N+\beta}{N-2\alpha+\beta}$.

\end{proposition}

\begin{proposition}\label{pr5}  \cite[Proposition 2.3]{CFV}
Let $\alpha\in(0,1]$ and $\beta\in [0, \alpha]$, then the mapping $f\mapsto \mathbb G_\alpha[f]$ is compact from  $L^{1}(\Omega,\rho^\beta_{\partial\Omega} dx)$ into $L^{q}(\Omega)$ for any $q\in [1,\frac{N}{N+\beta-2\alpha})$.
Moreover, for  $q\in [1,\frac{N}{N+\beta-2\alpha})$, there exists $c_{45}>0$ such that for any $f\in L^{1}(\Omega,\rho^{\beta}_{\partial\Omega}dx)$
  \begin{equation}\label{power1}
  \norm{\mathbb G_\alpha[f]}_{L^q(\Omega)}\leq c_{45}\norm f_{L^{1}(\Omega,\rho^{\beta}_{\partial\Omega}dx)}.
\end{equation}
 \end{proposition}

For $\nu\in\mathfrak{M}^b_{\partial\Omega}(\bar\Omega)$, $\nu_t$ is given in section 2.2 for $t\in(0,\sigma_0)$.  Let
$t_j=\frac1j\in(0,\sigma_0/4)$ if $j\ge j_0$ for some $j_0>0$.
Choose $\{\tilde\nu_n\}_n\subset C_0^1(\Omega)$  a sequence of nonnegative functions such that
supp$(\tilde\nu_n)\subset\Omega_{t_{j_0}-2^{-n}}\setminus\Omega_{t_{j_0}+2^{-n}}$ and
$\tilde\nu_{n}\to\nu_{t_{j_0}} $ in the  duality sense with $C(\bar
\Omega)$. Denote
$$\nu_{n,j}(x)=\left\{
\arraycolsep=1pt
\begin{array}{lll}
\tilde \nu_{n}(x+t_j\vec{n_x}),\quad& {\rm if}\quad x\in \Omega_{t_{j_0}-2^{-n}}\setminus\Omega_{t_{j_0}+2^{-n}}, \\[2mm]
0,&{\rm if\ not.}
\end{array}
\right.
$$

\begin{lemma}\label{lm 6.1}
Up to subsequence, we have that
$\nu_{n,j_n}\to\nu $ in  the duality sense with $C(\bar
\Omega)$, that is,
\begin{equation}\label{06-08}
  \lim_{n\to\infty}\int_{\bar \Omega}\zeta \nu_{n,j_n }dx=\int_{\bar \Omega}\zeta d\nu,\qquad\forall \zeta\in C(\bar \Omega).
\end{equation} Moreover,
$${\rm supp}(\nu_n)\subset \Omega_{\frac{t_n}2}\setminus \Omega_{2t_n}.$$
\end{lemma}
{\bf Proof.}
For any fixed $j$ and $\zeta\in C(\bar\Omega)$,
we observe that
\begin{eqnarray*}
\lim_{n\to\infty}\int_{\bar \Omega}\zeta \nu_{n,j}dx = \int_{\Omega}\zeta d\nu_{t_j}
\end{eqnarray*}
and pass $j\to\infty$, we derive that
\begin{eqnarray*}
\lim_{j\to\infty}\lim_{n\to\infty}\int_{\bar \Omega}\zeta \nu_{n,j}dx = \int_{\Omega}\zeta d\nu.
\end{eqnarray*}
The second argument is obvious by the definition of $\nu_{n,j}$.
\qquad$\Box$

\subsection{ Sub-linear }
In this subsection, we are devoted to prove the existence of weak solution to (\ref{eq 1.1})
when the source nonlinearity is sub-linear.
\smallskip

\noindent{\bf Proof of Theorem \ref{teo 3}  $(i)$.}
 Let $\{\nu_n\}$ be a sequence of nonnegative functions such that
$\nu_{n}\to\nu $ in  sense of duality with $C(\bar\Omega)$, see Lemma \ref{lm 6.1}.
By the Banach-Steinhaus Theorem, we may assume that $\norm{\nu_n}_{L^1(\Omega)}\le \norm{\nu}_{\mathfrak M^b (\Omega)}=1$ for all $n$.  We consider a sequence $\{g_n\}$ of $C^1$ nonnegative  functions defined on $\R_+$
such that $g_n(0)=g(0)$,
\begin{equation}\label{06-08-1}
  g_n\le g_{n+1}\le g,\quad \sup_{s\in\R_+}g_n(s)=n\quad{\rm and}\quad \lim_{n\to\infty}\norm{g_n-g}_{L^\infty_{loc}(\R_+)}=0.
\end{equation}
We set
$$M(v)=\norm{v}_{L^{1}(\Omega)}.$$

{\it Step 1. To prove that for $n\geq 1$,
 \begin{equation}\label{002.3}
 \arraycolsep=1pt
\begin{array}{lll}
 (-\Delta)^\alpha u= g_{n}(u)+kt_n^{-\alpha}\nu_n\quad & {\rm in}\quad\Omega,\\[2mm]
 \phantom{   (-\Delta)^\alpha }
u=0\quad & {\rm in}\quad \Omega^c
\end{array}
 \end{equation}
admits a nonnegative solution $u_n$ such that
$$M(u_n)\le \bar\lambda,$$
where $\bar\lambda>0$ independent of $n$. }

To this end, we define the operators $\{\mathcal{T}_n\}$ by
 $$\mathcal{T}_nu=\mathbb{G}_\alpha\left[g_n(u)+k t_n^{-\alpha}\nu_n\right],\qquad \forall u\in L^1_+(\Omega),$$
 where $L^1_+(\Omega)$ is the positive cone of $L^1(\Omega)$.
By  (\ref{power1})  and (\ref{06-08-2}), we have that
\begin{equation}\label{23-05-0}
 \arraycolsep=1pt
\begin{array}{lll}
  M(\mathcal{T}_nu)\le c_{45}\norm{g_n(u)+k t_n^{-\alpha}\nu_n}_{L^1 (\Omega,\rho_{\partial\Omega}^{\alpha}dx)}
  \\[2mm] \phantom{---- }
  \le c_3c_{45} \int_{\Omega}u^{p_0}\rho^\alpha(x)dx+c_{46}(k+\epsilon)
  \\[2mm] \phantom{---- }
  \le c_3c_{47}\int_{\Omega}u^{p_0} dx+c_{46}(k+\epsilon)
  \\[2mm] \phantom{---- }
  \le c_3 c_{48}(\int_{\Omega}u dx)^{p_0}+c_{46}(k+\epsilon)
  \\[2mm] \phantom{---- }
=c_3c_{48}M(u)^{p_0}+c_{46}(k+\epsilon),
\end{array}
\end{equation}
where $c_{47},c_{48}>0$ independent of $n$.
Therefore, we derive that
$$
  M(\mathcal{T}_nu)\le c_3c_{48} M(u)^{p_0}+c_{45}(k+\epsilon).
$$

If we assume that $M(u)\le \lambda$ for some $\lambda>0$, it implies
$$
 M(\mathcal{T}_nu)\le c_3c_{48} \lambda^{p_0}+c_{45}(k+\epsilon).
$$
 In the case of $p_0<1$,   the equation
$$
c_3c_{48}\lambda^{p_0 }+c_{45}(k+\epsilon)=\lambda
$$
admits a unique positive root $\bar\lambda$.
In the case of $p_0=1$,  for $c_3>0$ satisfying $c_3c_{48}<1$, the equation
$$
c_3c_{48}\lambda+c_{45}(k+\epsilon)=\lambda
$$
admits a unique positive root $\bar\lambda$.
For $M(u)\le \bar\lambda$, we obtain that
  \begin{equation}\label{07-05-5jingxuan}
  M(\mathcal{T}_nu)\le c_3c_{48}\bar\lambda^{p_0}+c_{45}(k+\epsilon)= \bar\lambda.
  \end{equation}
Thus, $\mathcal{T}_n$ maps $L^1(\Omega)$ into itself. Clearly, if $u_m\to u$ in $L^1(\Omega)$ as $m\to\infty$, then $g_n(u_m)\to g_n(u)$ in $L^1(\Omega)$ as $m\to\infty$, thus $\mathcal{T}_n$ is continuous.
For any fixed $n\in\N$, $\mathcal{T}_nu_m=\mathbb{G}_\alpha\left[g_n(u_m)+k \nu_n\right]$ and $\{g_n(u_m)+k \nu_n\}_m$ is
uniformly bounded in $L^1(\Omega,\rho_{\partial\Omega}^\beta dx)$, then it follows by Proposition \ref{pr5} that
$\{\mathbb{G}_\alpha\left[g_n(u_m)+k t_n^{-\alpha}\nu_n\right]\}_m$ is pre-compact in $L^1(\Omega)$, which implies that
$\mathcal{T}_n$ is a compact operator.

Let
$$
\displaystyle\begin{array}{lll}\displaystyle
\mathcal{G}=\{u\in L^1_+(\Omega): \ M(u)\le \bar\lambda \},
\end{array}
$$
  which is a closed and convex
set of $L^1(\Omega)$.  It infers by (\ref{07-05-5jingxuan}) that
$$\mathcal{T}_n(\mathcal{G})\subset \mathcal{G}.$$
 It follows by Schauder's fixed point theorem that there exists some $u_n\in L^1_+(\Omega)$ such that
$\mathcal{T}_nu_n=u_n$ and $M(u_n)\le \bar\lambda,$
where $\bar\lambda>0$ independent of $n$.

We observe that $u_n$ is a classical solution of  (\ref{002.3}). Let  open set $O$ satisfy $ O\subset \bar O\subset \Omega$.
By  \cite[Proposition 2.3]{RS}, for $\theta\in(0,2\alpha)$, there exists $c_{49}>0$ such that
$$\norm{u_n}_{C^{\theta}(O)}\le c_{49}\{\norm{g(u_n)}_{L^\infty(\Omega)}+kt_n^{-\alpha}\norm{\nu_n}_{L^{\infty}(\Omega)}\},$$
then applied  \cite[Corollary 2.4]{RS}, $u_n$ is $C^{2\alpha+\epsilon_0}$ locally in $\Omega$ for some $\epsilon_0>0$.
Then $u_n$ is a classical solution of (\ref{002.3}).
Moreover, from \cite[Lemma 2.2]{CV2}, we derive that
\begin{equation}\label{5.60000}
\int_\Omega u_n(-\Delta)^\alpha\xi dx=\int_\Omega g(u_n)\xi dx+k\int_\Omega\xi
t_n^{-\alpha}\nu_ndx,\quad \forall\xi\in \mathbb{X}_{\alpha}.
\end{equation}


{\it Step 2. Convergence. } We observe that  $\{g_n( u_n)\}$ is uniformly bounded in $L^1(\Omega,\rho_{\partial\Omega}^\alpha dx)$, so is $\{\nu_n\}$.
By Proposition \ref{pr5}, there exist a subsequence $\{u_{n_k}\}$ and $u$ such that
$u_{n_k}\to u$ a.e. in $\Omega$ and in $L^1(\Omega)$, then by (\ref{06-08-2}), we derive that
$g_{n_k}(u_{n_k}) \to g( u)$  in $L^1(\Omega)$. Pass the limit of (\ref{5.60000}) as $n_k\to \infty$ to derive that
 $$\int_\Omega u(-\Delta)^\alpha\xi=\int_\Omega g(u)\xi dx+k\int_\Omega\frac{\partial^\alpha\xi}{\partial \vec{n}^\alpha} d\nu,\quad \forall \xi\in\mathbb{X}_\alpha, $$
thus $u$ is a weak solution of (\ref{eq 1.1}) and $u$ is nonnegative since $\{u_n\}$ are nonnegative.
\qquad$\Box$

\subsection{Integral subcritical }
In this subsection, we prove   the existence of weak solution to (\ref{eq 1.1})
when the nonlinearity is integral subcritical.

\smallskip

\noindent{\bf Proof of Theorem \ref{teo 3}  $(ii)$.}
Let  $\{\nu_n\}\subset C^1(\bar \Omega)$ be a sequence of nonnegative functions given as the above
 and  $\norm{\nu_n}_{L^1(\Omega)}\le 2\norm{\nu}_{\mathfrak M ^b(\bar\Omega)}=1$ for all $n$.
We consider a sequence $\{g_n\}$ of $C^1$ nonnegative  functions defined on $\R_+$
satisfying $g_n(0)=g(0)$ and (\ref{06-08-1}).
We set
$$M_1(v)=\norm{v}_{M^{\frac{N+\alpha}{N-\alpha}}(\Omega,\rho^\alpha_{\partial\Omega} dx)}\quad{\rm and}\quad M_2(v)=\norm{v}_{L^{p_*}(\Omega)},$$
where $p_*$ is (\ref{1.4}).
We may assume that $p_*\in(1, \frac{N}{N-\alpha})$.
In fact, if $p_*\ge \frac{N}{N-\alpha}$, then  for any given $p\in(1, \frac{N}{N-\alpha})$, (\ref{1.4}) implies that
$$g(s)\le c_{4}s^p+\epsilon,\quad \forall s\in[0,1].$$

{\it Step 1. To prove that for $n\geq 1$,
 \begin{equation}\label{2.3}
 \arraycolsep=1pt
\begin{array}{lll}
 (-\Delta)^\alpha u= g_{n}(u)+kt_n^{-\alpha}\nu_n\quad & {\rm in}\quad\Omega,\\[2mm]
 \phantom{   (-\Delta)^\alpha_n }
u=0\quad & {\rm in}\quad \Omega^c
\end{array}
 \end{equation}
admits a nonnegative solution $u_n$ such that
$$M_1(u_n)+M_2(u_n)\le \bar\lambda,$$
where $\bar\lambda>0$ independent of $n$. }

To this end, we define the operators $\{\mathcal{T}_n\}$ by
 $$\mathcal{T}_nu=\mathbb{G}_\alpha\left[g_n(u)+kt_n^{-\alpha} \nu_n\right],\qquad \forall u\in L^1_+(\Omega).$$
By  Proposition \ref{general}, we have
\begin{eqnarray}
  M_1(\mathcal{T}_nu) &\le& c_{44}\norm{g_n(u)+kt_n^{-\alpha} \nu_n}_{L^1 (\Omega,\rho^{\alpha}_{\partial\Omega}dx)}\nonumber\\[2.5mm]
   &\le & c_{44} [\norm{g_n(u)}_{L^1(\Omega,\rho^{\alpha}_{\partial\Omega}dx)}+k].\label{06-08-10}
\end{eqnarray}
In order to deal with $\norm{g_n(u)}_{L^1(\Omega,\rho^{\beta}_{\partial\Omega}dx)}$, for $\lambda
>0$ we set
$$S_\lambda=\{x\in\Omega:u(x)>\lambda\}\quad {\rm and}\quad
\omega(\lambda)=\int_{S_\lambda}\rho^\alpha_{\partial\Omega} dx,$$
\begin{equation}\label{chenyuhang1}
\displaystyle\begin{array}{lll}
\displaystyle\norm{g_n(u)}_{L^1(\Omega,\rho^{\alpha}_{\partial\Omega}dx)}\le \int_{S^c_1}g(u)\rho^{\alpha}_{\partial\Omega}dx+\int_{S_1}
g(u)\rho^{\alpha}_{\partial\Omega}dx.
\end{array}
\end{equation}
We first deal with $\int_{S_1} g(u)\rho^{\alpha}dx$.  In fact, we observe that
$$\int_{S_1} g(u)\rho^{\alpha}_{\partial\Omega}dx=\omega(1) g(1)+\int_1^\infty \omega(s)dg(s),$$
where
$$\int_1^\infty  g(s)d\omega(s)=\lim_{T\to\infty}\int_1^T   g(s)d\omega(s).
$$
It infers by  Proposition \ref{pr 1} and Proposition \ref{general} that
there exists  $c_{50}>0$ such that
\begin{equation}\label{2.4}
\omega(s)\leq c_{50}M_1(u)^{\frac{N+\alpha}{N-\alpha}}s^{-\frac{N+\alpha}{N-\alpha}}
\end{equation}
and by  (\ref{1.4}) and Lemma \ref{lm 08-09} with $p=\frac{N+\alpha}{N-\alpha}$, there exist a sequence of increasing numbers $\{T_j\}$ such that $T_1>1$ and
$T_j^{-\frac{N+\alpha}{N-\alpha}}  g(T_j)\to 0$ when $j\to\infty$, thus
$$\displaystyle\begin{array}{lll}
\displaystyle \omega(1)  g(1)+ \int_1^{T_j} \omega(s)d g(s) \le c_{50}M_1(u)^{\frac{N+\alpha}{N-\alpha}} g(1)+c_{50}M(u)^{\frac{N+\alpha}{N-\alpha}}\int_1^{T_j} s^{-\frac{N+\alpha}{N-\alpha}}d g(s)
\\[4mm]\phantom{-----}\displaystyle
\leq c_{50}M_1(u)^{\frac{N+\alpha}{N-\alpha}}{T_j}^{-\frac{N+\alpha}{N-\alpha}}
g(T_j)+\frac{c_{50}M_1(u)^{\frac{N+\alpha}{N-\alpha}}}{\frac{N+\alpha}{N-\alpha}+1}\int_1^{T_j}
s^{-1-\frac{N+\alpha}{N-\alpha}} g(s)ds.
\end{array}$$
Therefore,
\begin{equation}\label{06-08-11}
\displaystyle\begin{array}{lll}
\int_{S_1} g(u)\rho^{\alpha}dx=\omega(1)g(1)+ \int_1^\infty \omega(s)\ dg(s)
\\[3mm]\phantom{------}
\leq \frac{c_{50}M_1(u)^{\frac{N+\alpha}{N-\alpha}}}{\frac{N+\alpha}{N-\alpha}+1}\int_1^\infty s^{-1-\frac{N+\alpha}{N-\alpha}} g(s)ds
\\[3mm]\phantom{------}
\displaystyle \le c_{50}g_\infty M_1(u)^{\frac{N+\alpha}{N-\alpha}},
\end{array}
\end{equation}
where $c_{50}>0$ independent of $n$.

We next deal with $  \int_{S^c_1}g(u)\rho^{\alpha}_{\partial\Omega}dx$.
For $p_*\in(1, \frac{N}{N-2\alpha+\beta})$, we have that
\begin{equation}\label{4.1}
\displaystyle\begin{array}{lll}
 \int_{S^c_1}g(u)\rho^{\alpha}_{\partial\Omega}dx\le c_{4}\int_{S_1^c}u^{p_*}\rho^{\alpha}_{\partial\Omega}dx+\epsilon\int_{S_1^c}\rho^{\alpha}_{\partial\Omega}dx
 \\[3mm]\phantom{------}
 \le c_{4}c_{51}\int_{\Omega}u^{p_*}dx+c_{51}\epsilon
  \\[3mm]\phantom{------}
 \leq c_{4}c_{51}M_2(u)^{p_*} +c_{51}\epsilon,
\end{array}
\end{equation}
where $c_{51}>0$ independent of $n$.

Along with (\ref{06-08-10}), (\ref{chenyuhang1}), (\ref{06-08-11}) and (\ref{4.1}), we derive
\begin{equation}\label{05-09-4}
  M_1(\mathcal{T}_nu)\le c_{44}c_{50}g_\infty M_1(u)^{\frac{N+\alpha}{N-\alpha}}+c_{44}c_{4}c_{51}M_2(u)^{p_*}+c_{44}c_{51}\epsilon+c_{44}k.
\end{equation}
By \cite[Theorem 6.5]{NPV} and (\ref{power1}), we derive that
$$
\displaystyle\begin{array}{lll}
M_2(\mathcal{T}_nu)\le c_{45}\norm{g_n(u)+k \nu_n}_{L^1 (\Omega,\rho_{\partial\Omega}^{\alpha}dx)},
\end{array}
$$
which along with  (\ref{chenyuhang1}), (\ref{06-08-11}) and (\ref{4.1}), implies that
\begin{equation}\label{4.3}
  M_2(\mathcal{T}_nu)\le c_{45}c_{50}g_\infty M_1(u)^{\frac{N+\alpha}{N-\alpha}}+c_{45}c_{4}c_{51}M_2(u)^{p_*}+c_{45}c_{51}\epsilon+c_{45}k.
\end{equation}

Therefore,  inequality (\ref{05-09-4}) and (\ref{4.3}) imply that
$$
M_1(\mathcal{T}_nu)+M_2(\mathcal{T}_nu)\le c_{52}g_\infty M_1(u)^{\frac{N+\alpha}{N-\alpha}}+c_{53}c_4M_2(u)^{p_*}+c_{54}\epsilon+c_{54}k,
$$
where $c_{52}=(c_{44}+c_{45})c_{50}$, $c_{21}=(c_{44}+c_{45})c_{51}$ and $c_{54}=c_{44}+c_{45}$.
If we assume that $M_1(u)+M_2(u)\le \lambda$, implies
$$
 M_1(\mathcal{T}_nu)+M_2(\mathcal{T}_nu)\le c_{52}g_\infty\lambda^{\frac{N+\alpha}{N-\alpha}}+c_{21}\lambda^{p_*}+c_{21}\epsilon+c_{54}k.
$$
Since $\frac{N+\alpha}{N-\alpha},\ p_*>1$, then there exist $k_0>0$ and $\epsilon_0>0$ such that for any $k\in(0,k_0]$ and $\epsilon\in(0,\epsilon_0]$,  the equation
$$
c_{52}g_\infty\lambda^{\frac{N+\alpha}{N-\alpha}}+c_{21}\lambda^{p_*}+c_{21}c_3\epsilon+c_{54}k=\lambda
$$
admits the largest root $\bar\lambda>0$.

We redefine $M(u)=M_1(u)+M_2(u)$, then for  $M(u)\le \bar\lambda$, we obtain that
 \begin{equation}\label{2.2}
  M(\mathcal{T}_nu)\le c_{52}g_\infty\bar\lambda^{\frac{N+\alpha}{N-\alpha}}+c_{21}\bar\lambda^{p_*}+c_{21}\epsilon+c_{54}k= \bar\lambda.
 \end{equation}
Especially, we have that
$$\norm{\mathcal{T}_nu}_{L^1(\Omega)}\le c_8M_1(\mathcal{T}_nu)|\Omega|^{\frac{2\alpha}{N+\alpha}}\le c_{23} \bar\lambda\quad{\rm if}\quad M(u)\le \bar\lambda.$$
Thus, $\mathcal{T}_n$ maps $L^1(\Omega)$ into itself. Clearly, if $u_m\to u$ in $L^1(\Omega)$ as $m\to\infty$, then $g_n(u_m)\to g_n(u)$ in $L^1(\Omega)$ as $m\to\infty$, thus $\mathcal{T}_n$ is continuous.
For any fixed $n\in\N$, $\mathcal{T}_nu_m=\mathbb{G}_\alpha\left[g_n(u_m)+k \nu_n\right]$ and $\{g_n(u_m)+k \nu_n\}_m$ is
uniformly bounded in $L^1(\Omega,\rho^\alpha dx)$, then it follows by Proposition \ref{pr5} that
$\{\mathbb{G}_\alpha\left[g_n(u_m)+k \nu_n\right]\}_m$ is pre-compact in $L^1(\Omega)$, which implies that
$\mathcal{T}_n$ is a compact operator.

Let
$$
\displaystyle\begin{array}{lll}\displaystyle
\mathcal{G}=\{u\in L^1_+(\Omega): \ M(u)\le \bar\lambda \}
\end{array}
$$
  which is a closed and convex
set of $L^1(\Omega)$.  It infers by (\ref{2.2}) that
$$\mathcal{T}_n(\mathcal{G})\subset \mathcal{G}.$$
 It follows by Schauder's fixed point theorem that there exists some $u_n\in L^1_+(\Omega)$ such that
$\mathcal{T}_nu_n=u_n$ and $M(u_n)\le \bar\lambda,$
where $\bar\lambda>0$ independent of $n$.

In fact, $u_n$ is a classical solution of  (\ref{2.3}).   Let  $O$ an open set satisfying $ O\subset \bar O\subset \Omega$.
By  \cite[Proposition 2.3]{RS}, for $\theta\in(0,2\alpha)$, there exists $c_{55}>0$ such that
$$\norm{u_n}_{C^{\theta}(O)}\le c_{55}\{\norm{g(u_n)}_{L^\infty(\Omega)}+kt_n^{-\alpha}\norm{\nu_n}_{L^{\infty}(\Omega)}\},$$
then applied  \cite[Corollary 2.4]{RS}, $u_n$ is $C^{2\alpha+\epsilon_0}$ locally in $\Omega$ for some $\epsilon_0>0$.
Then $u_n$ is a classical solution of (\ref{2.3}).
Moreover,
\begin{equation}\label{5.6}
\int_\Omega u_n(-\Delta)^\alpha\xi dx=\int_\Omega g(u_n)\xi dx+k\int_\Omega\xi
\nu_ndx,\quad \forall\xi\in \mathbb{X}_{\alpha}.
\end{equation}

{\it Step 2. Convergence. } Since  $\{g_n( u_n)\}$ and $\{\nu_n\}$ are uniformly bounded in $L^1(\Omega,\rho^\beta_{\partial\Omega} dx)$,
then by  Propostion \ref{pr5},  there exist a subsequence $\{u_{n_k}\}$ and $u$ such that
$u_{n_k}\to u$ a.e. in $\Omega$ and in $L^1(\Omega)$, and
$g_{n_k}(u_{n_k}) \to g( u)$ a.e. in $\Omega$.

 Finally we prove that  $g_{n_k}( u_{n_k})\to g( u)$ in $L^1(\Omega,\rho^\beta_{\partial\Omega} dx)$.
For $\lambda
>0$, we set $S_\lambda=\{x\in\Omega:|u_{n_k}(x)|>\lambda\}$  and
$\omega(\lambda)=\int_{S_\lambda}\rho^{\alpha}_{\partial\Omega}dx$, then for any Borel
set $E\subset\Omega$, we have that
\begin{equation}\label{chenyuhang1000}
\displaystyle\begin{array}{lll}
\displaystyle\int_{E}|g_{n_k}(u_{n_k})|\rho^\beta_{\partial\Omega} dx=\int_{E\cap
S^c_{\lambda}}g(u_{n_k})\rho^\beta_{\partial\Omega} dx+\int_{E\cap S_{\lambda}}
g(u_{n_k})\rho^\beta_{\partial\Omega} dx\\[4mm]\phantom{\int_{E}|g(u_{n_k})|\rho^{\alpha_{\partial\Omega}}dx}
\displaystyle\leq \tilde g(\lambda)\int_E\rho^\beta_{\partial\Omega} dx+\int_{S_{\lambda}} g(u_{n_k})\rho^\beta_{\partial\Omega} dx\\[4mm]\phantom{\int_{E}g(u_{n_k})\rho^\beta_{\partial\Omega} dx}
\displaystyle\leq \tilde
g(\lambda)\int_E\rho^\beta_{\partial\Omega} dx+\omega(\lambda) g(\lambda)+\int_{\lambda}^\infty \omega(s)d
g(s),
\end{array}
\end{equation}
where $\tilde g(\lambda)=\max_{s\in[0,\lambda]}g(s)$.

On the other hand,
$$\int_{\lambda}^\infty g(s)d\omega(s)=\lim_{T_m\to\infty}\int_{\lambda}^{T_m} g(s)d\omega(s).
$$
where $\{T_m\}$ is a sequence increasing number such that
$T_m^{-\frac{N+\alpha}{N-\alpha}} g(T_m)\to 0$ as $m\to\infty$,
which could obtained by assumption (\ref{1.4}) and Lemma \ref{lm 08-09} with $p={\frac{N+\alpha}{N-\alpha}}$.

It infers by  (\ref{2.4}) that
$$\displaystyle\begin{array}{lll}
\displaystyle \omega(\lambda)  g(\lambda)+ \int_{\lambda}^{T_m} \omega(s)d g(s) \le c_{50}  g(\lambda)\lambda^{-\frac{N+\alpha}{N-\alpha}}+c_{56}\int_{\lambda}^{T_m} s^{-\frac{N+\alpha}{N-\alpha}}d  g(s)
\\[4mm]\phantom{-----\ \int_{\lambda}^{T_m}  g(s)d\omega(s)}\displaystyle
\leq c_{56}T_m^{-\frac{N+\alpha}{N-\alpha}}g(T_m)+\frac{c_{56}}{\frac{N+\alpha}{N-\alpha}+1}\int_{\lambda}^{T_m}
s^{-1-\frac{N+\alpha}{N-\alpha}} g(s)ds,
\end{array}$$
where $c_{56}=c_{50}\frac{N+\alpha}{N-\alpha}$.
Pass the limit of $m\to\infty$, we have that
$$\omega(\lambda)  g(\lambda)+ \int_{\lambda}^\infty \omega(s)\ d  g(s)\leq \frac{c_{56}}{\frac{N+\alpha}{N-\alpha}+1}\int_{\lambda}^\infty s^{-1-\frac{N+\alpha}{N-\alpha}}  g(s)ds.
$$

Notice that the above quantity on the right-hand side tends to $0$
when $\lambda\to\infty$. The conclusion follows: for any
$\epsilon>0$ there exists $\lambda>0$ such that
$$\frac{c_{56}}{\frac{N+\alpha}{N-\alpha}+1}\int_{\lambda}^\infty s^{-1-\frac{N+\alpha}{N-\alpha}} g(s)ds\leq \frac{\epsilon}{2}.
$$
Since $\lambda$ is fixed, together with (\ref{chenyuhang1}), there exists $\delta>0$ such that
$$\int_E\rho^{\alpha}_{\partial\Omega}dx\leq \delta\Longrightarrow  g(\lambda)\int_E\rho^{\alpha}_{\partial\Omega}dx\leq\frac{\epsilon}{2}.
$$
This proves that $\{g\circ u_{n_k}\}$ is uniformly integrable in
$L^1(\Omega,\rho^\beta_{\partial\Omega} dx)$. Then $g\circ u_{n_k}\to g\circ u$ in
$L^1(\Omega,\rho^\beta_{\partial\Omega} dx)$ by Vitali convergence theorem.

Pass the limit of (\ref{5.6}) as $n_k\to \infty$ to derive that
 $$\int_\Omega u(-\Delta)^\alpha\xi=\int_\Omega g(u)\xi dx+k\int_\Omega\frac{\partial^\alpha\xi}{\partial \vec{n}^\alpha} d\nu,\quad \forall \xi\in\mathbb{X}_\alpha, $$
thus $u$ is a weak solution of (\ref{eq 1.1}) and $u$ is nonnegative since $\{u_n\}$ are nonnegative.
\qquad$\Box$

\bigskip
\bigskip

Huyuan Chen
\medskip

\noindent Department of Mathematics, Jiangxi Normal University,

\noindent Nanchang, Jiangxi 330022, PR China

  and

\noindent  Institute of Mathematical Sciences,  New York University  Shanghai,

\noindent Shanghai 200120, PR China

 \bigskip

 Hichem Hajaiej

 \medskip
\noindent  Institute of Mathematical Sciences,  New York University  Shanghai,

 \noindent Shanghai 200120, PR China

\end{document}